%% file: main.tex
\documentclass[a4paper,12pt]{article}

\usepackage{amsmath,
            mleftright,
            comment,
            amssymb,
            amsthm,
            nicefrac,
            xcolor,
            enumerate,
            authblk,
            color,
            etoolbox,
            mathrsfs,
            mathtools,
            bbm,
            xparse,
            stmaryrd,
            enumitem,
            mathtools,
            marginnote,
            xurl,
            graphicx,
            scalerel,
            relsize,
            stackengine,
            caption
            }

\usepackage{float}

\usepackage{tikz}
\usetikzlibrary{positioning, shapes.geometric, calc, fit, backgrounds}

\usepackage[colorlinks=true]{hyperref}

\usepackage{comment}
\usepackage[sort,capitalize]{cleveref}

\crefformat{equation}{(#2#1#3)}
\crefname{enumi}{item}{items}
\crefname{equation}{}{}

\usepackage[left=20mm,right=18mm,top=1.5cm,bottom=1cm,includeheadfoot]{geometry}

\usepackage[T1]{fontenc}
\usepackage[utf8]{inputenc}
\usepackage[english]{babel}
\usepackage[sort,capitalize]{cleveref}

\usepackage{csquotes}

\newtheorem{theo}{Theorem}[section]

\newtheorem{lemma}[theo]{Lemma}
\newtheorem{assumption}[theo]{Assumption}
\newtheorem{prop}[theo]{Proposition}
\newtheorem{corollary}[theo]{Corollary}
\newtheorem{setting}[theo]{Setting}

\theoremstyle{definition}
\newtheorem{ex}[theo]{Example}

\newtheorem{defi}[theo]{Definition}

\newtheorem{remark}[theo]{Remark}

\numberwithin{equation}{section}

\definecolor{darkgreen}{rgb}{0,0.6,0}

\newcommand{\Rn}{\mathbb{R}^n}
\newcommand{\Rq}{\mathbb{R}^q}
\newcommand{\Rm}{\mathbb{R}^m}
\newcommand{\net}{\mathrm{NN}}

\newcommand{\numberSteps}{N}

\newenvironment{mproof}[1]{\noindent \textit{Proof of {#1}.}}{\hfill \qed}

\title{
Convexity and strict convexity for compositional neural networks in high-dimensional optimal control
}
\author{
    Lars Gr\"une$^{1}$, Konrad Kleinberg$^{2}$, Thomas Kruse$^{3}$, Mario Sperl$^{4}$ \bigskip \\
    \small{$^1$ Mathematical Institute, University of Bayreuth,}
    \vspace{-0.1cm}\\
    \small{Germany; e-mail: \texttt{lars.gruene}\textcircled{\texttt{a}}\texttt{uni-bayreuth.de}}\smallskip\\
    \small{$^2$ Department of Mathematics \& Informatics, University of Wuppertal,}
    \vspace{-0.1cm}\\
    \small{Germany; e-mail: \texttt{kleinberg}\textcircled{\texttt{a}}\texttt{uni-wuppertal.de}}\smallskip\\
    \small{$^3$ Department of Mathematics \& Informatics, University of Wuppertal,}
    \vspace{-0.1cm}\\
    \small{Germany; e-mail: \texttt{tkruse}\textcircled{\texttt{a}}\texttt{uni-wuppertal.de}} \smallskip\\
    \small{$^4$ Mathematical Institute, University of Bayreuth,}
    \vspace{-0.1cm}\\
    \small{Germany; e-mail: \texttt{mario.sperl}\textcircled{\texttt{a}}\texttt{uni-bayreuth.de}}
}
\date{\today}

\input{commands.tex}

\begin{document}

\maketitle

\begin{abstract}
    Neural networks (NNs) have emerged as powerful tools for solving high-dimensional optimal control problems. 
In particular, their compositional structure has been shown to enable efficient approximation of high-dimensional functions, helping to mitigate the curse of dimensionality in optimal control problems.
In this work, we build upon the theoretical framework developed by Kang \& Gong (SIAM J.\ Control Optim. 60(2):786-813, 2022), particularly their results on NN approximations for compositional functions in optimal control. 
Theorem 6.2 in Kang \& Gong (SIAM J.\ Control Optim. 60(2):786-813, 2022) establishes that, under suitable assumptions on the compositional structure and its associated features, optimal control problems with strictly convex cost functionals admit a curse-of-dimensionality-free approximation of the optimal control by NNs. We extend this result in two directions. First, we analyze the strict convexity requirement on the cost functional and demonstrate that reformulating a discrete-time optimal control problem with linear transitions and stage costs as a terminal cost problem ensures the necessary strict convexity. Second, we establish a generalization of Theorem 6.2 in Kang \& Gong (SIAM J.\ Control Optim. 60(2):786-813, 2022) which provides weak error bounds for optimal control approximations by NNs when the cost functional is only convex rather than strictly convex. 

\end{abstract}

\newpage

\tableofcontents

\section{Introduction}

High-dimensional problems arise naturally in a multitude of applications, e.g., partial differential equations, (stochastic) optimal control, and stability of dynamical systems.
For such tasks, many conventional numerical algorithms are subject to the curse of dimensionality (COD) (see, e.g., \cite{Bellman57,NovakRitter1997,NovakWozniakowski2008I}). This means that their computational complexity to obtain a solution grows exponentially with the problem dimension and thus renders such algorithms impractical in high dimensions.
The utilization of machine learning, in particular through (deep) neural networks (NNs), has demonstrated empirical success in various high-dimensional tasks which were previously intractable due to the COD (see, e.g., \cite{han2018solving,ChanMikaelWarin2019,EHanJentzen2017CMStat,BeckEJentzen2017,beck2021Kolmogorov,mnih2015human,silver2016mastering,lillicrap2015continuous}).

The COD is particularly pronounced in optimal control when grid or mesh-based methods are employed to solve the associated Hamilton–Jacobi–Bellman (HJB) equations.
Early contributions already investigated the use of neural networks to approximate optimal value functions or feedback laws in nonlinear control problems (see, e.g., \cite{sontag1991feedback,narendra1990identification,chen1994adaptively}). More recently, this line of research has been extended to high-dimensional settings, where neural networks have been employed to approximate optimal value functions and compute feedback laws in problems with state dimensions beyond the reach of traditional grid-based approaches (see, e.g., \cite{kunisch2021semiglobal,kunisch2023optimal, nakamura2022,nakamura2021adaptive,sperl2025separable}). 

It is important to emphasize that, although neural networks (NNs) are universal approximators for continuous functions (see, e.g., \cite{Cybenko1989,hornik1989multilayer}), in general, NN approximations of functions with high-dimensional input are subject to the COD (see, e.g., \cite{Mhaskar1996NNoptimal,Pinkus1999MLP,Guehring2023Expressivity}).
Thus, the objective of attaining theoretical guarantees for when and why NNs overcome the COD has inspired a broad and active field of research. A variety of such results have been obtained in different contexts: for instance, for (control) Lyapunov functions \cite{gruene2021_ComputingLyapunovFunctions,sperl2025_CLF}, for partial (integro-)differential equations \cite{hutzenthaler2019proof,gonon2023deep}, and for problems in (stochastic) optimal stopping and control \cite{ReisingerZhang2020zerosum,kang2022feedforward,Gonon2024DeepStopping,jentzen2025DNNBellman}.

One particular mechanism that allows NNs to mitigate the COD is when the target function exhibits favorable structural properties, such as a compositional structure in terms of lower-dimensional components. The role of compositionality in enabling efficient approximation has been emphasized in \cite{poggio2017_WhyAndWhenCanDeepAvoidCoD} and investigated further in various directions (see, e.g., \cite{kang2022feedforward,kang2024data,gong2023approximationReLU,riekert2023deep,dahmen2025compositional_published,gruene2021_ComputingLyapunovFunctions,sperl2025_CLF}).

In this paper, we explore the role of compositionality in optimal control problems to obtain curse-of-dimensionality-free representations using neural networks. To this end, we build on the theory of neural network approximation for compositional functions developed in \cite{kang2022feedforward}, with a particular focus on the results in Section 6 in this reference. It has been established in \cite[Theorem 6.2]{kang2022feedforward} that, under suitable assumptions on the compositional structure and its so-called compositional features, optimal control problems with terminal costs (also known as Mayer term) that are strictly convex in the control admit a curse-of-dimensionality-free representation of the optimal control.
We extend this result in two main directions. First, we examine the strict convexity assumption on the terminal cost. We demonstrate that the reformulation of a linear optimal control problem with stage costs as a terminal cost problem ensures the necessary strict convexity. Second, we prove a generalization of \cite[Theorem 6.2]{kang2022feedforward} that requires only convexity in the control input, rather than strict convexity. This relaxation comes at the expense that the network yields a control that only approximates the optimal value (a so-called ``weak approximation'') rather than a control that approximates the optimal trajectory (a so-called ``strong approximation'').

The outline of the paper is as follows: In Section \ref{sec:preliminaries} we collect some preliminary results from the literature.
In \cref{sec:sep_lin_disc}, we present a neural network approximation result for discrete time optimal control problems with separated stage costs. This can be seen as a distillation of a problem class for which \cite[Theorem~6.2]{kang2022feedforward} can be applied.
In Section \ref{sec:NN_convex_problems} we prove the extensions of Theorem 6.1 and Theorem 6.2 in \cite{kang2022feedforward} for cost functionals that are convex with respect to the control input, respectively.
Finally, Section \ref{sec:conclusion} concludes the paper.

%%%%%%%%%%%%%%%

We close this introductory section by establishing notation which is employed throughout the remainder of the paper.
For a function $\psi \colon \mathbb{R}^n \times \mathbb{R}^q \to \mathbb{R}^q$, $(x,u) \mapsto \psi(x,u)$, and for $k \in \mathbb{N}$, we denote for $x \in \mathbb{R}^n$ and $u \in \mathbb{R}^q$ the iterates $\psi^k(x,u) \coloneqq \psi(x, \psi^{k-1}(x,u))$ with the convention $\psi^0(x,u) \coloneqq u$. Further, for a function $f \colon \R^n \times \R^q \to \R^n$ and $U = (u_1, \dots, u_N) \in \R^{qN}$ for some $N \in \N$, we define for any  $k \in \{1, \dots, N\}$ the $k$-fold composition of $f$ via $f^{k+1}(x,\cdot)(U) := f\big(f^{k}(x,\cdot)(U),\, u_k\big)$ with $f^{0}(x,\cdot)(U) = x$. Unless explicitly stated otherwise, for any vector $x \in \mathbb{R}^n$ and any matrix $A \in \mathbb{R}^{n \times n}$ we write $\|x\|$ and $\|A\|$ for the Euclidean norm of $x$ and the corresponding induced norm of $A$, respectively. For $x \in \mathbb{R}^n$ and $r \geq 0$ we define $\mathcal{B}_r(x) \coloneqq \{ y \in \mathbb{R}^n \mid \|x - y \| \leq r \}$, i.e., the closed Euclidean ball of radius $r$ centered at $x$. Further, for $n,m \in \mathbb{N}$ and a compact set $\dom \subset \mathbb{R}^n$ let $W^\infty_{m} \pr{\dom}$ denote the space of all functions $f \colon \Omega \to \R$ whose derivatives up to order $m$ exist almost everywhere and satisfy for every multiindex $\alpha \in \N^n$ with $\abs{\alpha} \le m$ that $\norm{D_\alpha f}_{\infty, \dom} < \infty$, where $\norm{\cdot}_{\infty, \dom}$ denotes the essential supremum norm.

%%%%%%%%%%%%%%%

%%%%%%%%%%%%%%%%% 

\section{The compositional framework and applications in optimal control} \label{sec:preliminaries}

\subsection{Preliminaries on compositional functions}
Throughout this paper, we use the concept of compositional functions as introduced in \cite[Definition 3.1]{kang2022feedforward} and discussed thereafter. To ensure that our exposition is self-contained, we provide a brief summary of the main aspects in the following.  

\begin{defi}[Definition 3.1 in \cite{kang2022feedforward}] \label{def:compositional_function}
    A compositional function is a triplet $\pr{f, \cG^f, \cL^f}$, where
    \begin{enumerate}[label=(\roman*)]
        \item $\cG^f = (\cV^f, \cE^f)$ is a directed acyclic graph with node set $\cV^f$ and edge set $\cE^f \subset \cV^f \times \cV^f$. 
        \item $\cL^f \colon \cV^f \to \N$ assigns a layer number for each node in $\cV^f$ such that 
        \[
        (g,h) \in \mathcal{E}^f \quad \Rightarrow \quad \cL^f(g) < \cL^f(h).
        \]
        Nodes with no inward edge are called input nodes and have layer number $0$, while nodes with no outward edge are called output nodes and have layer number
        \[
        l^f_{\max} = \max_{g \in \cV^f} \cL^f(g).
        \]
        All remaining nodes are called hidden nodes. 
        \item For $0 \leq i \leq l^f_{\max}$ we denote the set of all nodes with layer number $i$ as $\{ f_{i,1}, \dots, f_{i, N_i} \}$. 
        Each node $f_{i,j}$ with $i > 0$, $ 1 \leq j \leq N_i$, is a function
         \[
           f_{i,j} \colon \mathbb{R}^{d_{i,j}} \to \mathbb{R},
        \]
       where $d_{i,j}$ is the number of incoming edges to $f_{i,j}$ in $\cG^f$. 
       \item With $d := N_0$, we assign an input $x \in \R^d$ to the input nodes $f_{0,1}, \dots f_{0,d}$. Then, recursively we define for all remaining nodes $f_{i,j}$
        \[
        z_{i,j}(x) = f_{i,j}\big( (z_{i-1,k}(x))_{(f_{i-1,k},f_{i,j}) \in \cE^f} \big),
        \]
        where the vector $(z_{i-1,k}(x))_{(f_{i-1,k},f_{i,j}) \in \cE^f}$ collects the outputs of all predecessors $f_{i-1,k}$ of $f_{i,j}$.
        \item Finally, with $q := N_{l^f_{\max}}$, the function $f \colon \R^d \to \R^q$ is given by 
        \[
        f(x) = \big(z_{f_{l^f_{\max}, 1}}(x), \dots, z_{f_{l^f_{\max}, q}}(x)\big).
        \]
    \end{enumerate}
\end{defi}

\begin{figure}[H]
    \centering
    \includegraphics[scale=0.8]{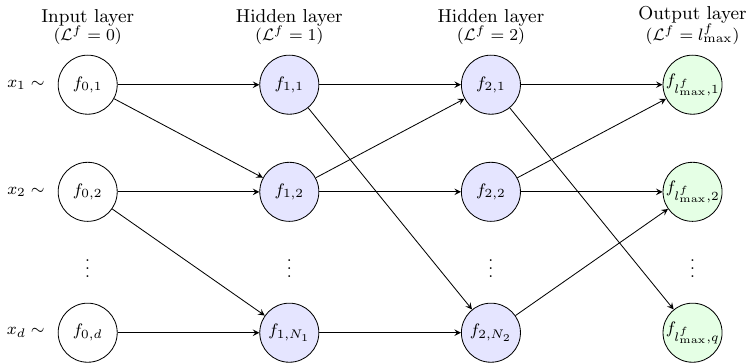}
    \caption{Illustration of a compositional function.}
    \label{fig:def_compositional_function}
\end{figure}

Figure~\ref{fig:def_compositional_function} illustrates the structure of a compositional function as introduced in Definition~\ref{def:compositional_function}. For simplicity, we occasionally refer to a compositional function $\pr{f, \cG^f, \cL^f}$ simply as $f$ in the following. For a compositional function, we denote the set of nodes that are either an input node or that represent a linear affine function by $\cV_L^f$. The set of general nodes is defined as $\cV_G^f := \cV^f \setminus \cV_L^f$. Altogether, this notion of compositional functions allows to express how a mapping $f$ is built from the lower dimensional functions $f_{i,j}$ of its nodes in $\cG^f$. Throughout this paper, we adopt the following assumption on compositional functions, which is taken from \cite[Assumption 1]{kang2022feedforward}. 

\begin{assumption} \label{asm:compositional_functions}
  We assume that $\pr{f, \cG^f, \cL^f}$ is a compositional function, where $f$ maps from $[ -R, R]^d$ to $ \mathbb{R}^q $ for some $ R > 0 $ and $ d, q \in \N_{> 0}$. Further, we assume that the nodes $f_{i,j} \in \cV^f_G $ are functions in $ W^\infty_{m_{i,j}}([-R_{i,j}, R_{i,j}]^{d_{i,j}}) $ for some $ m_{i,j}, d_{i,j} \in \N_{>0} $ and $ R_{i,j} > 0 $. The ranges and domains of all nodes are compatible for composition, i.e., if $ (f_{i,j}, f_{l,k}) $ is an edge in $ \cG^f $, then the range of $ f_{i,j} $ is contained in the interior of the domain of $ f_{l,k} $.
\end{assumption}

In the spirit of Definition~4.3 in \cite{kang2022feedforward}, we view a neural network ($\net$) as a compositional function of a particular type. 

\begin{defi}[Definition~4.3 in \cite{kang2022feedforward}] \label{def:neural_network}
    A neural network with activation function $\sigma$ is a compositional function $\pr{f^\net, \cG^{f^\net}, \cL^{f^\net}}$, where each function $f_{i,j}^\net \in \cV_G^{f^\net}$ is of the form 
    \begin{equation*}
        f_{i,j}^\net \colon \R^{d_{i,j}} \to \R, \quad x \mapsto \sigma( w^T x + b), 
    \end{equation*}
    where $w \in \R^{d_{i,j}}$ and $b \in \R$ are parameters assigned to $f_{i,j}^\net$ and $x\in\R^{d_i,j}$ is the vector of input values to the node $f_{i,j}^\net$. Further, each function in the output layer is given as a weighted sum of the nodes in the previous layer. The complexity of a $\net$ is defined as total number of hidden nodes in $\cG^{f^\net}$ denoted by $\size\pr{f^\net} \in \N$. 
\end{defi}

A neural network with one hidden layer is called a shallow neural network. For compositional functions satisfying Assumption \ref{asm:compositional_functions} a quantitative approximation result with neural networks has been shown in \cite[Theorem 4.10]{kang2022feedforward}. There, the main idea is to replace each general node of the given compositional function with a shallow neural network. This results in a neural network that approximates the given compositional function. The resulting approximation error depends on the so-called features of a compositional function. For the convenience of the reader, we re-state the definition of the compositional features and the approximation result here: 

\begin{defi}[Definition 4.1 in \cite{kang2022feedforward}] \label{def:comp_features}
    Let $\pr{f, \cG^f, \cL^f}$ satisfy Assumption \ref{asm:compositional_functions}. Then the numbers $\pr[\big]{ r_{\max}^f, \Lambda^f, L_{\max}^f, \abs{\cV_G^f} } \in (0,\infty)^4$ are called \textit{features} of $\pr{f, \cG^f, \cL^f}$ if  $r_{\max}^f, \Lambda^f, L_{\max}^f \in \R_{>0}$ satisfy for all $f_{i,j} \in \cV_G^f$ that 
    \begin{align*}
        d_{i,j} m_{i,j}^{-1} & \le r_{\max}^f, \\ 
        \max \cu*{R_{i,j}^{m_{i,j}}, 1} \sum_{| \alpha | \le m_{i,j}} \| \diff_\alpha f_{i,j} \|_{\infty, \dom_{i,j}} & \le \Lambda^f, \\ 
        \abs{L_{i,j}^f} & \le L_{\max}^f,  
    \end{align*}
    where $\dom_{i,j} = [- R_{i,j},R_{i,j}]^{d_{i,j}}$ and $L_{i,j}^f$ is a Lipschitz constant associated to $f_{i,j}$, see \cite[Definition 3.6]{kang2022feedforward}. 
\end{defi}
In the remainder of this paper, the notation $\pr[\big]{ r_{\max}^f, \Lambda^f, L_{\max}^f, \abs{\mathcal{V}_{\mathcal{G}}^f} }$ will always indicate compositional features of a compositional function $\pr{f, \mathcal{G}^f, \mathcal{L}^f}$, even if this is not explicitly mentioned.

\begin{theo}[Theorem 4.10 in \cite{kang2022feedforward}] \label{thm:NN_approximation_compositional}
    Let $\pr{f, \cG^f, \cL^f}$ satisfy Assumption \ref{asm:compositional_functions} and let $\pr[\big]{ r_{\max}^f, \Lambda^f, L_{\max}^f, \abs{\cV_G^f}}$ be features of $\pr{f, \cG^f, \cL^f}$. Then for any $n_w \in \N$ there exists a $\net$ $f^\net$ with $\size\pr{f^\net} = n_w \abs{\cV_G^f}$ such that for all $x \in [-R,R]^d$ it holds that 
    \begin{equation*}
        \norm{ f(x) - f^\net(x) } \le C \Lambda^f L_{\max}^f \abs{\cV_G^f} n_w^{- \frac{1}{r_{\max}^f}}, 
    \end{equation*}
    where $C > 0$ is a constant determined by $\cu{ \pr{ d_{i,j}, m_{i, j} } \in \N^2 : f_{i,j} \in \cV_G^f}$.
\end{theo}

\subsection{Optimal control with compositionality}
We aim to leverage the fact that compositional functions can be approximated well with neural networks for high-dimensional optimal control problems. Our starting point is the result stated in Theorem 6.2 in \cite{kang2022feedforward}. There, the authors consider an optimal control problem with dynamics $\dot{x} = f(x,u)$ over a fixed time horizon $t \in [0,T]$ together with a terminal cost functional of the form $J(x,U) = \Psi(\Phi(x,U))$, where $\Psi \colon \Rn \to \R_{\ge 0}$ and $\Phi(x,U)$ denotes the state of the system starting at $x$ and and controlled using a piecewise constant control over $N \in \N$ time steps specified by a vector $U \in \mathbb{R}^{qN}$. The statement of \cite[Theorem 6.2]{kang2022feedforward}  is that, under suitable assumptions, the optimal control can be approximated with a $\net$ whose complexity grows sub-exponentially with the dimension. The main ingredients required for the result are:
\begin{enumerate}
    \item The compositional structure of both $f$ and $\Psi$,
    \item Suitable assumptions on the growth of the involved compositional features, 
    \item Positive-definiteness of the Hessian of $J(x,U)$ with respect to $U$ at all considered initial values $x$. \label{asm:hessian_J}
\end{enumerate}

We refer to Section 6 in \cite{kang2022feedforward} for the detailed description of the theorem and its assumptions.
The following example shows that Item \ref{asm:hessian_J} above, i.e., strict convexity of $J(x,U)$ with respect to $U$ does not hold in general, even if the dynamics are linear and $\Psi$ is strictly convex.  

\begin{ex}\label{ex:linear_dynamic_strict_convex_terminal_cost}
    Let $a, b \in \R$,
    let $\Psi \colon \R \to \R$ satisfy for all $x \in \R$ that $\Psi^{''}(x) > 0$, let $\numberSteps \in \N$.
    We consider a discrete-time optimal control problem with $\numberSteps$ time steps,  dynamics governed by
    \begin{equation}
        x^+ = ax + bu, \qquad x \in \R, \quad u \in \R,
    \end{equation}
    and terminal costs $\Psi$, allowing for discrete-time control inputs $U = \br{u_0, \dots, u_{\numberSteps-1}}^\top \in \R^\numberSteps$, where the control $u_k \in \R$ is applied at stage $k$. By induction one shows that the terminal stage of the system given an initial state $x \in \R$ and a sequence of control inputs $U = \br{u_0, \dots, u_{\numberSteps-1}}^\top \in \R^\numberSteps$ can be described by $\Phi \colon \R \times \R^\numberSteps \to \R$, $\Phi(x, U) = a^\numberSteps x + b \sum_{j = 0}^{\numberSteps - 1} a^{\numberSteps - 1 - j}u_j$. The associated cost functional to be mimimized is thus given by $J(x, U) = \Psi ( \Phi(x, U))$. The Hessian of $J$ w.r.t. $U \in \R^\numberSteps$ is of the form
    \begin{equation}
        \diff^2_U J(x,U) = v v^\top b^2 \Psi^{''}(\Phi(x, U)), \qquad x \in \R, \quad U \in \R^\numberSteps,
    \end{equation}
    where $v = \pr{a^{N - 1}, a^{N - 2}, \dots, 1}^\top \in \R^\numberSteps$.
    Note that for all $\numberSteps \ge 2$, $x \in \R$, $U \in \R^\numberSteps$ it holds that $\det(\diff^2_U J(x,U)) = 0$, since the Hessian is a rank-1 matrix. Hence it holds for all $\numberSteps \ge 2$, $x \in \R$, $U \in \R^\numberSteps$ that $\diff_U^2 J (x, U) \not\succ 0$. Moreover, it holds for all $\numberSteps \ge 2$, $x \in \R$, $U \in \R^\numberSteps$, $w \in \R^\numberSteps$ that $w^\top \diff_U^2 J(x, U) w \ge 0$. Thus, $J$ is convex but not strictly convex w.r.t. $U$.
\end{ex}

\cref{ex:linear_dynamic_strict_convex_terminal_cost} motivates the investigation of the two following directions: First, in \cref{sec:sep_lin_disc}, we introduce stage costs to gain strict convexity in the cost functional for linear dynamics. Second, we study how the approximation result can be extended to settings where only convexity, but not strict convexity, is present. In \cref{sec:NN_convex_problems}, we weaken the strict convexity requirement of Theorem~6.2 in \cite{kang2022feedforward} and establish an approximation result for convex setups. Since optimal controls are not necessarily unique without strict convexity, the result is stated in a weak form: instead of comparing the neural network output with \emph{the} optimal control, we compare its performance with the performance of \emph{an} optimal control, where performance is measured via the cost functional.

\begin{remark}
    While our focus in this work is on discrete-time optimal control problems, we note that the extension to continuous-time dynamics is also possible. In particular, Theorem 5.1 in \cite{kang2022feedforward} shows that continuous-time systems can be approximated via Euler discretization in a manner that preserves the compositional structure relevant for neural network approximation.
\end{remark}

\section{Enforcing strict convexity by introducing stage costs}
\label{sec:sep_lin_disc}

In this section, we consider the following discrete-time optimal control problem with stage costs.

\begin{setting} \label{set:general_dynamic}
    Let $n, q \in \N$, let $f \colon \R^n \times \R^q \to \R^n$ a function which describes the system dynamic by 
    \begin{equation}\label{eq:ext_state_space_dynamic}
        x^+ = f\pr{x,u}, \qquad x \in \R^n, \quad u \in \R^q,
    \end{equation}
    let $\numberSteps \in \N$ be the time horizon,
    for $x \in \R^n$, $U \in \R^{q\numberSteps}$ let the sequence $\pr{x_{U, k}}_{k \in \cu{0,1,\dots, \numberSteps}}$ denote the state evolution with initial state $x \in \R^n$ and control input\footnote{Zero-order hold (ZOH) control input (see, e.g.,\ \cite{kang2022feedforward}).} $U =\begin{bmatrix} u_0^\top & u_1^\top & \cdots & u_{\numberSteps - 1}^\top \end{bmatrix} ^\top \in \R^{q\numberSteps}$,
    let $\cost \colon \Rn \times \Rq \to \R$ denote stage costs,
    let $g \colon \R^n \to \R$ denote terminal costs,
    let the cost functional $J \colon \R^n \times \R^{q\numberSteps} \to \R$ satisfy for all $x \in \R^n$, $U = \begin{bmatrix} u_0^\top & u_1^\top & \cdots & u_{\numberSteps - 1}^\top \end{bmatrix} ^\top \in \R^{q\numberSteps}$ that 
    \begin{equation} \label{eq:ext_state_space_stage_cost_problem}
        J(x,U) = g\pr{ x_{U, \numberSteps} } + \sum_{k = 0}^{\numberSteps - 1} \cost(x_{U,k}, u_k).
    \end{equation}
\end{setting}

\subsection{Strict convexity for problems with linear dynamics and separated stage costs}
\label{subsec:strict_convex_disc_lin}

The following remark provides a representation of the state in case of linear dynamics. This is utilized in proving strict convexity of the cost functional in the case of separated stage costs (see \cref{prop:sep_cost_strong_convexity} below).

\begin{remark} \label{rem:disc_lin_roll_off}
    Assume \cref{set:general_dynamic}, let $A \in \R^{n \times n}$, let $B \in \R^{n \times q}$,
    assume for all $x \in \R^n$, $u \in \R^q$ that $f(x, u) = Ax + Bu$,
    let $C^{\pr{k}}_j \in \R^{n \times q}$, $k \in \cu{0,1,\dots, \numberSteps}$, $j \in \cu{0, 1, \dots, \numberSteps - 1}$, satisfy for all $k \in \cu{0,1,\dots, \numberSteps}$, $j \in \cu{0, 1, \dots, \numberSteps - 1}$ that
    \begin{equation}
        \begin{split}
            C^{\pr{k}}_j = \begin{cases}
                A^{k - j - 1} B & \colon k \in \cu{1, 2,\dots,\numberSteps} \text{ and } j \in \cu{0,1\dots,k - 1} \\
                0               & \colon (k = 0) \vee (j \in \cu{k, \dots, \numberSteps - 1}),
            \end{cases}
        \end{split}
    \end{equation}
    let $\hat{C}^{\pr{k}} \in \R^{n \times (q\numberSteps)}$, $k \in \cu{0,1,\dots, \numberSteps}$, satisfy for all $k \in \cu{0,1,\dots, \numberSteps}$ that 
    \begin{equation}
        \hat{C}^{\pr{k}} = \begin{pmatrix} C^{\pr{k}}_0 & C^{\pr{k}}_1 & \cdots & C^{\pr{k}}_{\numberSteps - 1} \end{pmatrix}.
    \end{equation}
    Then it holds for all $x \in \R^n$, $U \in \R^{q\numberSteps}$, $k \in \cu{0,1,\dots,\numberSteps }$ that
    \begin{equation}
        x_{U, k} = A^kx + \SmallSum{j = 0}{k - 1}A^{k - 1 - j}Bu_j = A^k x + \SmallSum{j = 0}{k - 1} C^{\pr{k}}_j u_j = A^k x + \hat{C}^{\pr{k}} U.
    \end{equation}
\end{remark}

The next proposition establishes that, in the case of linear dynamics and separated stage costs, the associated cost functional is strictly convex with respect to the control input under the assumption that, both, the terminal costs and the state stage costs are convex, while the control stage costs are strictly convex.

\begin{prop}
    \label{prop:sep_cost_strong_convexity}
    Assume the situation as in \cref{rem:disc_lin_roll_off}, 
    assume that $\terminal$ is $\cC^2$,
    let $\cost_1\colon \Rn \to \R$, $\cost_2 \colon \Rq \to \R$ be $\cC^2$,
    assume that for all $x \in \R^n$, $u \in \R^q$ it holds that $\diff^2 \terminal(x) \succeq 0$, $\diff^2\cost_1(x) \succeq 0$, $\diff^2\cost_2(u) \succ 0$, and $\cost(x,u) = \cost_1(x) + \cost_2(u)$.
    Then $J$ is $\cC^2$ and it holds for all $x \in \R^n$, $U \in \R^{q\numberSteps}$ that $\diff^2_U J (x, U) \succ 0$.
\end{prop}

\begin{mproof}{\cref{prop:sep_cost_strong_convexity}}
    Let $\varphi_k \colon \R^n \times \R^{q\numberSteps} \to \R^n$, $k \in \cu{0,1,\dots, \numberSteps}$, satisfy for all $x \in \R^n$, $U \in \R^{q\numberSteps}$, $k \in \cu{0,1,\dots, \numberSteps}$ that $\varphi_k(x,U) = A^kx + \hat{C}^{\pr{k}} U$.
    First, observe that \cref{rem:disc_lin_roll_off} demonstrates for all $x \in \R^n$, $U \in \R^{q\numberSteps}$, $k \in \cu{0,1,\dots, \numberSteps}$ that $\varphi_k \pr{x, U} = x_{U, k}$.
    Next, note that the assumption that $g$, $\cost_1$, and $\cost_2$ are $\cC^2$
    and the fact that for all $k \in \cu{0, 1, \dots, \numberSteps}$ it holds that $\varphi_k$ is $\cC^2$ proves that $J$ is $\cC^2$.
    Throughout the remainder of this proof let $\bI_q \in \R^{q \times q}$, $q \in \N$, denote the $q$-dimensional identity matrix and let $p_k \colon \R^{q\numberSteps} \to \Rq$, $k \in \cu{0,1,\dots, \numberSteps - 1}$ satisfy for all $U = \begin{bmatrix}
        u_0^\top & u_1^\top & \cdots & u_{\numberSteps - 1}^\top
    \end{bmatrix}^\top \in \R^{q\numberSteps}$, $k \in \cu{0,1,\dots, \numberSteps - 1}$ that $p_k(U) = u_k$.
    Observe that for all $U \in \R^{q\numberSteps}$, $k \in \cu{0,1, \dots, \numberSteps - 1}$ it holds that $\diff p_k (U) = \begin{pmatrix}
        0 & \cdots & 0 & \bI_q & 0 & \cdots & 0
    \end{pmatrix} \in \R^{q \times q\numberSteps}$.
    Next, note that for all $k \in \cu{0,1,\dots, \numberSteps - 1}$, $U \in \R^{q\numberSteps}$ it holds that
    \begin{equation}
        \begin{split}
            \diff (\cost_2 \circ p_k) (U) &= \diff \cost_2 (p_k(U)) \diff p_k(U) = \diff \cost_2 (u_k)  \begin{pmatrix}
                0 & \cdots & 0 & \bI_q & 0 & \cdots & 0
            \end{pmatrix}
        \end{split}
    \end{equation}
    and $\diff^2 \pr{ \cost_2 \circ p_k } \pr{U} = \pr*{\begin{smallmatrix}
                0 & \cdots & 0 & \bI_q & 0 & \cdots & 0
            \end{smallmatrix}}^\top  \diff^2 \cost_2 (u_k) \pr*{ \begin{smallmatrix}
                0 & \cdots & 0 & \bI_q & 0 & \cdots & 0
            \end{smallmatrix}}$.
    This and the assumption that for all $u \in \Rq$ it holds that $\diff^2\cost_2 \pr{u} \succ 0$ establish for all $U \in \R^{q\numberSteps}$ that
    \begin{equation} \label{eq:sep_cost_strong_convexity_H_control_pos_definite}
        \begin{split}
            \sum_{k = 0}^{\numberSteps - 1} \diff^2 \pr{\cost_2 \circ p_k} \pr{U} = \pr*{\begin{smallmatrix}
                \diff^2 \cost_2 \pr{u_0}    & \phantom{0}   & 0                 \\
                \phantom{0}     & \ddots        & \phantom{0}       \\
                0               & \phantom{0}   & \diff^2 \cost_2 \pr{u_{\numberSteps - 1}}
            \end{smallmatrix}} \succ 0.
        \end{split}
    \end{equation}
    Moreover, note that for all $x \in \Rn$, $U \in \R^{q\numberSteps}$, $k \in \cu{0,1,\dots, \numberSteps - 1, \numberSteps}$ it holds that $\diff_U \varphi_k (x, U) = \hat{C}^{\pr{k}}$.
    This yields for all $x \in \R^n$, $U \in \R^{q\numberSteps}$, $k \in \cu{0,1,\dots, \numberSteps - 1}$ that
    \begin{equation}
        \begin{split}
            \diff_U\pr{ \cost_1 \circ \varphi_k } \pr{x, U} = \diff \cost_1 \pr{\varphi_k \pr{x, U}} \diff_U \varphi_k \pr{x, U} = \diff \cost_1 \pr{\varphi_k \pr{x, U}} \hat{C}^{\pr{k}}
        \end{split}
    \end{equation}
    and
    \begin{equation}
        \begin{split}
            \diff_U \pr{ \terminal \circ \varphi_\numberSteps} \pr{x, U} = \diff \terminal \pr{\varphi_\numberSteps \pr{x, U}} \diff \varphi_\numberSteps \pr{x, U} = \diff \terminal \pr{\varphi_\numberSteps \pr{x, U}} \hat{C}^{\pr{\numberSteps}}.
        \end{split}
    \end{equation}
    This shows for all $x \in \R^n$, $U \in \R^{q\numberSteps}$, $k\in \cu{0,1, \dots, \numberSteps - 1}$ that 
    \begin{equation}
        \begin{split}
            \diff^2_U \pr{\cost_1 \circ \varphi_k} \pr{x, U} = \pr{\hat{C}^{\pr{k}}}^\top \diff^2 \cost_1(\varphi_k(x, U)) \hat{C}^{\pr{k}} 
        \end{split}
    \end{equation}
    and
    \begin{equation}
        \begin{split}
            \diff^2_U \pr{\terminal \circ \varphi_\numberSteps} \pr{x, U} =  \pr{\hat{C}^{\pr{\numberSteps}}}^\top \diff^2 \terminal\pr{\varphi_\numberSteps(x, U)} \hat{C}^{\pr{\numberSteps}}.
        \end{split}
    \end{equation}
    This and the assumption that for all $x \in \R^n$ it holds that $\diff^2 \cost_1 \pr{x} \succeq 0$ and $\diff^2 \terminal \pr{x} \succeq 0$ ensures that for all $x \in \R^n$, $U \in \R^{q\numberSteps}$, $k \in \cu{0,1,\dots, \numberSteps - 1}$ it holds that $\diff^2_U \pr{\cost_1 \circ \varphi_k} \pr{x, U} \succeq 0$ and $\diff^2_U \pr{\terminal \circ \varphi_\numberSteps} \pr{x, U}\succeq 0$.
    Combining this with \eqref{eq:sep_cost_strong_convexity_H_control_pos_definite} implies for all $x \in \R^n$, $U \in \R^{q\numberSteps}$ that $\diff^2_U J\pr{x, U} \succ 0$. The proof of \cref{prop:sep_cost_strong_convexity} is thus complete.
\end{mproof}

\subsection{The extended state space problem and its compositional structure}
\label{subsec:ext_state_space_comp_struct}

In this section we reformulate the discrete-time optimal control problem from \cref{set:general_dynamic} as a problem involving only terminal costs and investigate the compositional structure of the resulting problem.
This structure then enables us to apply \cite[Theorem~6.2]{kang2022feedforward}, which we do in \cref{cor:stage_costs} in \cref{subsec:dnn_approx_stage_cost}.

We start by providing the reformulation, for which we introduce an extended state space.

\begin{remark} \label{rem:conv_reform_cost_funct}
    Assume \cref{set:general_dynamic}, let $\bff \colon \R^{n + 1} \times \R^q \to \R^{n + 1}$ and $\bfg \colon \R^{n + 1} \to \R$ be defined by
    \begin{equation}
        \begin{split}
            \bff \pr{x,y,u} = \begin{bmatrix}
                f(x,u) \\ y + \cost \pr{x,u}
            \end{bmatrix} \qquad \text{and} \qquad \bfg(x,y) = g(x) + y
        \end{split}
    \end{equation}
    for all $x \in \R^n$, $y \in \R$, $u \in \R^q$.
    We refer to $\bff$ as the extended system dynamics and to $\bfg$ as the extended terminal cost function.
    For $x \in \R^n$ and $U = \begin{bmatrix} u_0^\top & u_1^\top & \cdots & u_{\numberSteps - 1}^\top \end{bmatrix}^\top \in \R^{q\numberSteps}$ let the sequence $\pr{\overline{x}_{U, k}}_{k \in \cu{0,1, \dots, \numberSteps}}$ denote the state evolution in the extended state space with respect to the extended system dynamic $\bff$ with initial state $\pr{x, 0} \in \R^{n + 1}$ and control input $U = \begin{bmatrix} u_0^\top & u_1^\top & \cdots & u_{\numberSteps - 1}^\top \end{bmatrix}^\top \in \R^{q\numberSteps}$.
    Observe that for all $x \in \R^n$, $U = \begin{bmatrix} u_0^\top & u_1^\top & \cdots & u_{\numberSteps - 1}^\top \end{bmatrix}^\top \in \R^{q\numberSteps}$ it holds that
    \begin{equation}
        \overline{x}_{U, k} = \bff \pr{ \dots \bff \pr{ \bff\pr{x, 0, u_0}, u_1 }, \dots, u_{k - 1}} = \bff \pr{x, 0, \cdot}^{k} \pr{U} = \begin{bmatrix}
                x_{U,k} \\ \SmallSum{j = 0}{k - 1} \cost \pr*{x_{U,j}, u_j} 
                \end{bmatrix} 
    \end{equation}
    and $J\pr{x, U} = \bfg \pr{ \overline{x}_{U, \numberSteps}}$.
\end{remark}

The following lemma demonstrates that the reformulation provided by Remark \ref{rem:conv_reform_cost_funct} inherits the compositional structure of the defining functions. The key insight here is that the family of compositional functions is closed under certain algebraic operations such as parallelization, linear combination, and composition which is demonstrated in \cite[Section 3]{kang2022feedforward}. While the parallelization of compositional functions is not explicitly mentioned in \cite[Section 3]{kang2022feedforward}, it appears implicitly in the construction of linear combinations of compositional functions (cf.\ \cite[Definition~3.2]{kang2022feedforward}). Thus, we omit the proof of \cref{lem:comp_struct_disc_lin} and instead provide, in \cref{fig:extended_dynamic,fig:extended_terminal_cost} below, a visualization and an explanation of the constructed compositional functions appearing in \cref{lem:comp_struct_disc_lin}.

\begin{lemma} \label{lem:comp_struct_disc_lin}
Assume \cref{set:general_dynamic},
assume that $f$, $g$, and $\cost$ are compositional functions given by the triplets $\pr{f, \cG^f, \cL^f}$, $\pr{g, \cG^g, \cL^g}$, and $(\cost, \cG^\cost, \cL^\cost)$,
let $\bff \colon \R^{n + 1} \times \R^q \to \R^{n + 1}$ and $\bfg \colon \R^{n + 1} \to \R$ satisfy for all $x \in \R^n$, $y \in \R$, $u \in \R^q$ that 
\begin{equation}
    \begin{split} \label{eq:comp_struct_disc_lin}
        \bff\pr{x, y, u} = \begin{bmatrix}
            f\pr{x, u} \\ y + \cost(x, u)
        \end{bmatrix} \qquad \text{and} \qquad
            \bfg(x, y) = g(x) + y,
    \end{split}
\end{equation}
Then $\bff$ and $\bfg$ are compositional functions.
\end{lemma}

\bigskip \noindent

\begin{figure}[H]
    \centering
    \begin{minipage}[t]{0.45\textwidth}
        \centering
        \includegraphics[width=0.95\linewidth]{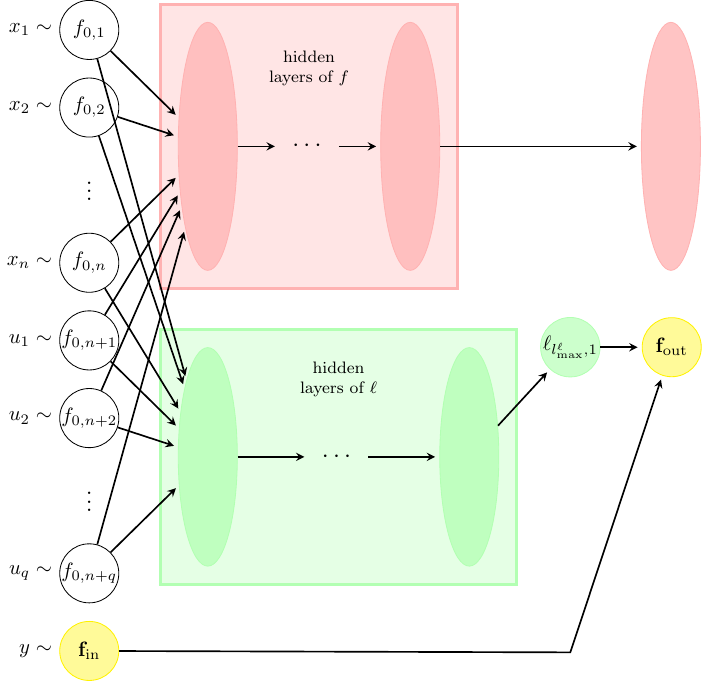}
        \caption{Extended dynamic $\bff$}
        \label{fig:extended_dynamic}
    \end{minipage}
    \begin{minipage}[t]{0.45\textwidth}
        \centering
        \includegraphics[width=0.95\linewidth]{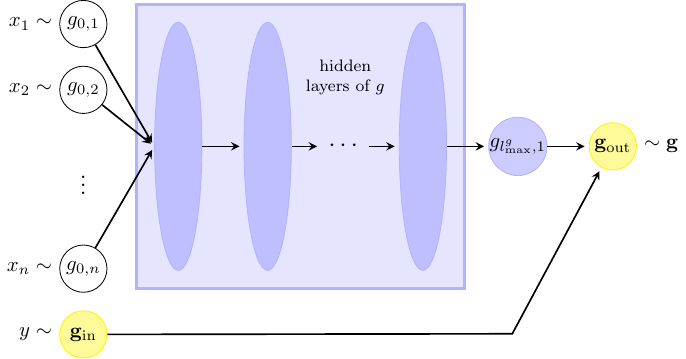}
        \caption{Extended terminal cost $\bfg$}
        \label{fig:extended_terminal_cost}
    \end{minipage}
\end{figure}

\cref{fig:extended_dynamic,fig:extended_terminal_cost} visualize the compositional graphs of the functions $\bff$ and $\bfg$ in \eqref{eq:comp_struct_disc_lin} in \cref{lem:comp_struct_disc_lin}. The extended system dynamics $\bff$ is constructed as the parallelization of the system dynamics $f$ and the stage cost function $\cost$ as in \cref{set:general_dynamic} together with the introduction of two new linear nodes corresponding to the additional input and output of the extended dynamic $\bff$. The compositional graph of the extended terminal cost function $\bfg$ is obtained by adding two linear nodes to the existing compositional graph of the terminal cost $g$ in \cref{set:general_dynamic} which correspond to the additional input of the terminal cost $\bfg$.

The next result, \cref{lem:comp_feat_disc_lin} below, derives compositional features of the extended system dynamics $\bff$ and the extended terminal cost function $\bfg$ of the reformulated system in \cref{rem:conv_reform_cost_funct,lem:comp_struct_disc_lin} based on the construction explained above and the compositional features of the dynamics $f$, stage cost function $\cost$, and terminal cost function $\terminal$ of the stage cost system as in \cref{set:general_dynamic}.

\begin{lemma} \label{lem:comp_feat_disc_lin}
    Assume the situation of \cref{lem:comp_struct_disc_lin},
    assume that the compositional functions $f, g$, and $\cost$ satisfy \cref{asm:compositional_functions}, 
    let $\pr*{ r_{\max}^h, \Lambda^h, L_{\max}^h, \cV_G^h } \in (0,\infty)^4$, $h \in \cu{f, g, \cost}$, satisfy for all $h \in \cu{f, g, \cost}$ that $\pr*{ r_{\max}^h, \Lambda^h, L_{\max}^h, \cV_G^h }$ are features of the compositional function $h$.
    Then 
    \begin{equation}
        \pr{\max \cu{r_{\max}^f, r_{\max}^\cost}, \max \cu{\Lambda^f, \Lambda^\cost}, \max \cu{ L_{\max}^f, L_{\max}^\cost}, \abs{\cV_G^f} + \abs{\cV_G^\cost}} \qquad \text{and} \qquad (r^g_{\max}, \Lambda^g, L_{\max}^g, \abs{\cV_G^g})
    \end{equation}
    are compositional features of $\bff$ and $\bfg$.
\end{lemma}

\begin{mproof}{\cref{lem:comp_feat_disc_lin}}
    Throughout this proof let $R^h \in (0,\infty)$, $h \in \cu{f, g, \cost}$,
    let $d^h \in \N$, $h \in \cu{f,g,\cost}$, satisfy $d^f = d^\cost = n + q$ and $d^g = n$,
    assume for all $h \in \cu{f, g, \cost}$ that the domain of $h$ is the hypercube $\br{-R^h, R^h}^{d^h}$,
    assume w.l.o.g.\ that $R^f = R^\cost$.
    Hence, restricting the range of the added input node in $\bff$ ($\bfg$ resp.) to the interval $\br{-R^f, R^f}$ ($\br{-R^g, R^g}$ resp.) ensures that the domain of $\bff$ ($\bfg$ resp.) is again a hypercube, specifically $\br{-R^f, R^f}^{d^f + 1}$ ($\br{-R^g, R^g}^{d^g + 1}$ resp.).
    Moreover, note that the added nodes in the construction of the compositional function $\pr{\bfg, \cG^\bfg, \cL^\bfg}$ are linear nodes.
    This ensures $\cV^g_G = \cV^\bfg_G$.
    Combining this
    with the fact that features of compositional functions depend only on their general nodes
    implies that $(r^g_{\max}, \Lambda^g, L_{\max}^g, \cV_G^g)$ are features of the compositional function $\bfg$.
    Next, observe that the added nodes
    in the construction of the compositional function $\pr{\bff, \cG^\bff, \cL^\bff}$ are linear nodes, i.e., it holds that $\cV_G^\bff = \cV_G^f \cup \cV_G^\cost$.
    Hence, it holds that $\pr{\max \cu{r_{\max}^f, r_{\max}^\cost}, \max \cu{\Lambda^f, \Lambda^\cost}, \max \cu{ L_{\max}^f, L_{\max}^\cost}, \abs{\cV_G^f} + \abs{\cV_G^\cost}}$
    are features of the compositional function $\pr{\bff, \cG^\bff, \cL^\bff}$.
    The proof of \cref{lem:comp_feat_disc_lin} is thus complete.
\end{mproof}

\subsection{Neural network approximations of optimal controls with stage costs}
\label{subsec:dnn_approx_stage_cost}
In this section, we present, in \cref{cor:stage_costs} below, a $\net$ approximation result for optimal control problems with linear dynamics and separated stage cost.
While in \cite[Theorem~6.2]{kang2022feedforward} it is assumed that the cost functional $J$, as the composition of the terminal cost and the dynamic, has a positive definite Hessian w.r.t.\ $U$, here we impose the convexity assumptions on the fundamental components of the optimal control problem, namely the stage and terminal cost functions, yielding a more modular set of assumptions.

\begin{corollary} \label{cor:stage_costs}
    Assume \cref{set:general_dynamic},
    let $A \in \R^{n \times n}$, $B \in \R^{n \times q}$,
    assume for all $x \in \R^n$, $u \in \R^q$ that $f(x,u) = Ax + Bu$,
    let $\cost_1 \colon \R^n \to \R$, $\cost_2 \colon \R^q \to \R$,
    assume that $\terminal$, $\cost_1$, and $\cost_2$ are compositional functions that satisfy \cref{asm:compositional_functions} and that all nodes in the compositional graphs of $\terminal$, $\cost_1$, and $\cost_2$ are $\cC^2$,
    assume for all $x \in \R^n$, $u \in \R^q$ that $\diff^2 \terminal\pr{x} \succeq 0$, $\diff^2\cost_1\pr{x} \succeq 0$, $\diff^2 \cost_2 \pr{u} \succ 0$, and $\cost(x,u) = \cost_1(x) + \cost_2(u)$,
    let $\dom \subseteq \R^n$ be compact,
    assume that for every $x \in \dom$ there exists $U^*(x) \in \R^{q\numberSteps}$ such that it holds that
    \begin{equation}
        \inf_{U \in \R^{q\numberSteps}} J(x,U) = J(x, U^*(x)),
    \end{equation}
    let $U_0 \in \R^{q\numberSteps}$,
    let $\gamma \in \pr{0, \infty}$,
    assume that for all $x \in \dom$ it holds that $U^*(x) \in \cB_{\gamma}(U_0)$.
    Then for every $\eps > 0$ there exists a $\net$ $U^{*\net} \colon \dom \to \R^{q\numberSteps}$ such that the inequality
    \begin{equation}
        \norm{U_\eps^{*\net} \pr{x} - U^* \pr{x} } \le 3 \eps
    \end{equation}
    holds for all $x \in \dom$ and the complexity of $U_\eps^{*\net}$ is bounded by
    \begin{equation} \label{eq:cor_stage_costs_estimate}
        \size \pr{ U_\eps^{*\net}} \le C \pr{q \pr{ \abs{\cV_G^{\cost_1}} + \abs{\cV_G^{\cost_2}} + \abs{\cV_G^{\terminal}} }}^{2r + 1} \eps^{- \pr{8r + 1}},
    \end{equation}
    for some constant $C > 0$ and $r = \max \cu{ 1, r_{\max}^{\cost_1}, r_{\max}^{\cost_2}, r_{\max}^{\terminal} }$.
    Here $r_{\max}^\terminal$, $r_{\max}^{\cost_1}$, and $r_{\max}^{\cost_2}$ are part of the compositional features of the functions $g \colon \br{-R, R}^n \to \R$, $\cost_1 \colon \br{-R, R}^n \to \R$, and $\cost_2 \colon \br{-R, R}^q \to \R$
    and $R > 1$ is a constant large enough such that $\dom \subseteq \br{-R, R}$, $\cB_{2\gamma}(U_0) \subseteq \br{- R, R}^{q \numberSteps}$ and for all $k \in \cu{1,\dots, \numberSteps}$, $x \in \dom$, $U \in \cB_{2\gamma}(U_0)$ it holds that
    \begin{equation}
        \begin{split}
            &x_{U, k} \in \br*{-\frac{R}{2}, \frac{R}{2}}^n \qquad \text{and} \qquad \pr*{ \sum_{j = 0}^{k - 1}  \cost_1(x_{U, j}) + \cost_2(u_j) } \in \br*{-\frac{R}{2}, \frac{R}{2}}.
        \end{split}
    \end{equation}
    The constant $C$ in the estimate \eqref{eq:cor_stage_costs_estimate} depends, in particular, on the compositional features of $g$, $\ell_1$, and $\ell_2$, and on the time horizon $N$, but is independent of the accuracy parameter $\eps$ and the dimension $n$.
\end{corollary}

\begin{mproof}{\cref{cor:stage_costs}}
    First, observe that $f \colon \R^n \times \R^q \ni (x, u) \mapsto Ax + Bu \in \R^n$ is a compositional function.
    All nodes of $f$ are linear nodes. 
    Its graph can be described without any hidden layers: The function $(x, u) \mapsto Ax + Bu$ as a linear function can be merged into a net with $n + q$ input nodes and $n$ output nodes and the node function associated to the $j$-th output node is given by $\R^n \times \R^q \ni (x, u) \mapsto A_{j, \cdot} x + B_{j, \cdot} u \in \R^n$ (cf.\ \cref{fig:compositional_structure_f}).
    \begin{center}
        \captionsetup{type=figure}
        \includegraphics[scale=0.8]{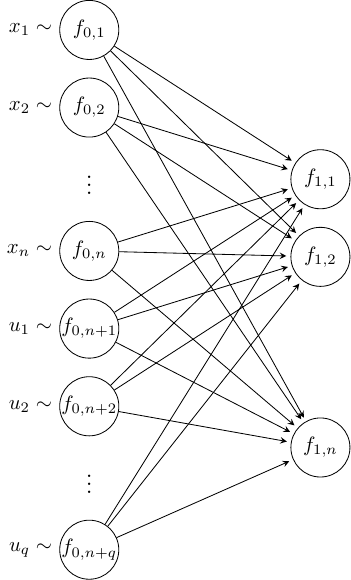}
        \captionof{figure}{Compositional structure of $f$}
        \label{fig:compositional_structure_f}
    \end{center}
    Thus all nodes of the compositional function $f$ are indeed linear nodes.
    Hence the set of general nodes is empty. 
    Consequently, one tuple of compositional features of $f$ (as it consists only of linear nodes and thus $\abs{\cV_G^f}$ = 0) can be chosen as $(r^f_{\max}, \Lambda^f, L^f_{\max}, \abs{ \cV_G^f }) = (1,0,0,0)$. 
    Next, observe that the assumption that $\cost_1$ and $\cost_2$ are compositional functions
    and the assumption that for all $(x, u) \in \R^n \times \R^q$ it holds that $\cost(x, u) = \cost_1(x) + \cost_2(u)$
    yield that $\cost$ is a compositional function and a tuple of features of $\cost$ is given by
    \begin{equation}
        \pr*{  \max\cu*{r_{\max}^{\cost_1}, r_{\max}^{\cost_2}},  \max \cu*{ \Lambda^{\cost_1}, \Lambda^{\cost_2} }, \max \cu*{ L_{\max}^{\cost_1}, L_{\max}^{\cost_2} }, \abs{ \cV_G^{\cost_1} } + \abs{\cV_G^{\cost_2}} }.
    \end{equation}
    Combining this,
    the compositional structure of $f$,
    \cref{lem:comp_struct_disc_lin},
    and \cref{lem:comp_feat_disc_lin} ensures that the function
    \begin{equation}
        \begin{split} \label{eq:stage_costs_extended_f}
            &\bff \colon \R^n \times \R \times \R^q \supseteq\br{-R, R}^{n + 1 + q} \to \R^{n + 1}\\
            &\bff \pr{x, y, u} = \begin{pmatrix}
                f(x, u) \\
                y + \cost(x, u)
            \end{pmatrix} = \begin{pmatrix}
                Ax + Bu \\
                y + \cost_1(x) + \cost_2(u)
            \end{pmatrix}, \quad (x,y,u) \in \br{-R, R}^{n + 1 + q} 
        \end{split}
    \end{equation}
    is a compositional function with a tuple of features given by
    \begin{equation}
        \begin{split}\label{eq:stage_costs_extended_f_features}
            \pr*{ \max\cu*{1, r_{\max}^{\cost_1}, r_{\max}^{\cost_2}},  \max \cu*{ \Lambda^{\cost_1}, \Lambda^{\cost_2} }, \max \cu*{ L_{\max}^{\cost_1}, L_{\max}^{\cost_2} }, \abs{ \cV_G^{\cost_1} } + \abs{\cV_G^{\cost_2}} }.
        \end{split}
    \end{equation}
    Similarly, it holds that the function
    \begin{equation}
        \begin{split}\label{eq:stage_costs_extended_g}
            \bfg\colon \R^n \times \R \supseteq \br{-R, R}^{n + 1} \to \R, \quad \bfg(x,y) = g(x) + y, \quad (x,y) \in  \br{-R, R}^{n + 1}
        \end{split}
    \end{equation}
    is also a compositional function with a tuple of features given by $\pr{r_{\max}^g, \Lambda^g, L_{\max}^g, \abs{\cV_G^g}}$, the features of $g$.
    Next,
    note that the assumption that for all $x \in \R^n$, $u \in \R^q$ it holds that $\diff^2 \terminal\pr{x} \succeq 0$, $\diff^2 \cost_1 \pr{x} \succeq 0$. $\diff^2 \cost_2 \pr{u} \succ 0$, and $\cost(x,u) = \cost_1(x) + \cost_2(u)$,
    and \cref{prop:sep_cost_strong_convexity}
    yield that for all $(x, U) \in \R^n \times \R^{q \numberSteps}$ it holds that $\diff^2_U J\pr{x, U} \succ 0$.
    Furthermore, 
    \cref{rem:conv_reform_cost_funct} ensures for all $x \in \dom$, $U \in \mathcal{B}_{2 \gamma} \pr{U_0}$ that $J\pr{x, U} = \bfg \pr{\overline{x}_{U, \numberSteps}}$.
    This and \cite[Theorem~6.2]{kang2022feedforward} yield that for every $\eps > 0$ there exists a deep neural network $U_\eps^{* \net} \colon \R^n \to \R^{q \numberSteps}$ such that for all $x \in \dom$ it holds that $\norm{ U_\eps^{*\net}\pr{x} - U^{*}\pr{x} } \le 3 \eps$ and the complexity of $U_\eps^{*\net}$ is bounded in the following sense:
    \begin{equation}
        \size \pr{ U_\eps^{*\net}} \le C \pr{q \pr{ \abs{\cV_G^{\cost_1}} + \abs{\cV_G^{\cost_2}} + \abs{\cV_G^{\terminal}} }}^{2r + 1} \eps^{- \pr{8r + 1}},
    \end{equation}
    for some constant $C > 0$. This concludes the proof of \cref{cor:stage_costs}.
\end{mproof}

\section{Weak approximations for non-strictly convex problems}
\label{sec:NN_convex_problems}

In this section, we address the second direction that emerged from our discussion following Example \ref{ex:linear_dynamic_strict_convex_terminal_cost}. Specifically, we study settings in which only convexity rather than strict convexity of the cost functional with respect to the control vector $U$ is present and prove a generalization of \cite[Theorem 6.2]{kang2022feedforward} to this case. As a result of this relaxed convexity assumption, the uniqueness of minimizers is no longer guaranteed. Consequently, the approximation of the optimal control by neural networks can only be established in a weak sense, that is, we measure the performance of the NN approximation through closeness of the cost functional values instead of comparing the resulting control and the optimal control as in the strictly convex setting. 

\subsection{Weak error estimates for gradient descent}

We begin by stating two auxiliary lemmata. The first one is a result for the convergence of the gradient descent method for convex functions.

\begin{lemma}[Theorem 2.1.14 and Corollary 2.1.15 in \cite{nesterov2018lectures}] \label{lem:gradient_descent}
    Let $f \colon \Rn \to \R$ be a convex $\cC^2$-function. Further, let the gradient $\nabla f \colon \Rn \to \Rn$ be Lipschitz-continuous with Lipschitz-constant $L>0$. Assume that there exists some $x^* \in \Rn$ such that 
    \begin{equation*}
        \inf_{x \in \Rn} f(x) = f(x^*). 
    \end{equation*}
    Define the mapping $\Psi \colon \Rn \to \Rn$, $x \mapsto x - \frac{1}{L} \nabla f(x)$. Then for all $k \in \N$ and $x \in \Rn$ it holds that 
            \begin{equation} \label{eq:gradient_descent_distance_x}
                \lVert \Psi^{k+1}(x) - x^* \rVert \leq \lVert \Psi^{k}(x) - x^* \rVert. 
            \end{equation}
        and 
            \begin{equation} \label{eq:gradient_descent_distance_f}
                f(\Psi^k(x)) - f(x^*) \leq \frac{2L}{k+4} \lVert x - x^* \rVert^2. 
            \end{equation}
\end{lemma}
Further, we need the following property of the Lipschitz constants of iterated maps. 
\begin{lemma} \label{lem:lipschitz_power_k}
    Let $\dom \subset \Rn$ and $f, \tilde{f} \colon \dom \to \dom$. Assume that $f$ is non-expansive, i.e.\ it is Lipschitz on $\dom$ with a Lipschitz constant $L \leq 1$. Suppose that there exists some $\varepsilon \in \R_{\geq 0}$ such that $\norm*{f(x) - \tilde{f}(x)} \leq \varepsilon$ for all $x \in \dom$. Then for all $k \in \N$ and $x \in \dom$ it holds that 
    \begin{equation*}
        \norm*{f^k(x) - \tilde{f}^k(x)} \leq k \varepsilon. 
    \end{equation*}
\end{lemma}
\begin{mproof}{\cref{lem:lipschitz_power_k}}
    The assertion holds for $k = 1$ by assumption. Suppose that it is true for $k-1$. Then we obtain for $x \in \dom$
    \begin{align*}
        \norm*{f^k(x) - \tilde{f}^k(x)} & \leq \norm*{f^k(x) - f(\tilde{f}^{k-1}(x))} + \norm*{f(\tilde{f}^{k-1}(x)) - \tilde{f}^k(x)} \\ & \leq \norm*{f^{k-1}(x) - \tilde{f}^{k-1}(x)} + \varepsilon \leq (k-1) \varepsilon + \varepsilon = k \varepsilon.  
    \end{align*}
\end{mproof} 

\subsection{Neural networks approximations for optimizers in non-strictly convex minimization problems}
In the following theorem, we address an optimization problem whose cost function is convex. Assuming that this function can be approximated arbitrarily well by a neural network, we construct a neural network that approximates an optimizer and establish an associated error bound together with the required network complexity. This result extends Theorem 6.1 in \cite{kang2022feedforward} to the case of convex functions and serves as the main building block for the subsequent result on neural network approximations of optimal controls in convex problems.

\begin{theo} \label{theo:approx_convex_cost}
    Let $\dom \subset \Rn$ be compact, $\cU \subset \Rm$ be convex, and $J \colon \dom \times \cU \to \R$, $(x,U) \mapsto J(x,U)$ be a $\cC^2$ function that is convex in $U$ for all $x \in \dom$. Assume that there exists some $U_0 \in \cU$ and $\gamma > 0$ with $\cB_{3\gamma}(U_0) \subset \cU$ such that for all $x \in \dom$ there exists some $U^*(x) \in \cB_\gamma(U_0)$ satisfying 
    \begin{equation*}
        \inf_{U \in \cU} J(x, U) = J(x, U^*(x)). 
    \end{equation*}
    Suppose that for each $\delta > 0$ there exists a NN $J_{\delta}^\net \colon \Rn \times \Rm \to \R$ with $ \size \pr{J^{\net}_\delta}$ neurons that approximates $J$ with accuracy $\delta > 0$ in the sense that for all $(x,U) \in \dom \times \cB_{2 \gamma}(U_0)$  
    \begin{equation} \label{eq:approx_convex_cost_NN_approx_J}
        \lvert J(x, U) - J_{\delta}^\net(x,U) \rvert \leq \delta. 
    \end{equation}
    Define the constants 
    \begin{equation} \label{eq:approx_convex_cost:L1_L2}
        L_1 \coloneq \max_{x \in \dom, \\ U \in \cB_{3\gamma}(U_0)} \norm*{\frac{\partial J}{\partial U}(x,U)}, \quad L_2 \coloneq \max_{x \in \dom, U \in \cB_{2\gamma}(U_0)} \norm*{\frac{\partial^2 J}{\partial U^2}(x,U)}. 
    \end{equation}
    Then for any $k \in \N_{>0}$ and $h, \delta \in \R_{>0}$ small enough such that  
    \begin{equation} \label{eq:approx_convex_cost_h_delta_small}
        k \pr*{h \sqrt{m} + \frac{2}{h L_2} \sqrt{m} \delta} \leq \gamma
    \end{equation}
    there exists a NN $U^{* \net}_{k,h,\delta} \colon \Rn \to \cB_{3 \gamma}(U_0)$ with $\size \pr{U^{*\net}_{k,h,\delta}} = 2 k m \size \pr{J^{\net}_\delta}$ neurons such that for all $x \in \dom$ we have 
    \begin{equation} \label{eq:approx_convex_cost_main_inequality}
         J(x,U^{*\net}_{k,h,\delta}(x)) - J(x, U^*(x)) \leq L_1 k \pr*{h \sqrt{m} + \frac{2}{h L_2} \sqrt{m} \delta} + \frac{2 L_2}{k+4} \gamma^2. 
    \end{equation}
\end{theo}
\begin{remark}
    Note that the NN $U^{*\net}_{k,h,\delta}$ in Theorem \ref{theo:approx_convex_cost} depends on the constants $k \in \N_{>0}$, $h \in \R_{>0}$, and $\delta \in \R_{>0}$. These parameters are crucial for two main reasons: On the one hand, the appropriate choice of these parameters allows to achieve an arbitrarily low approximation error in the right-hand side of equation \eqref{eq:approx_convex_cost_main_inequality}. On the other hand, in the case that $\cU \neq \Rm$ the choice of $k, h$, and $\delta$ is also essential to guarantee $U^{*\net}_{k,h,\delta}(\dom) \subset \cU$., i.e., that the range of the NN contains only admissible controls. 
\end{remark}

\begin{mproof}{\cref{theo:approx_convex_cost}}
    We want to utilize the fact that the gradient descent method is globally convergent for convex functions with a Lipschitz gradient, see \cref{lem:gradient_descent}. First note that the gradient $\frac{\partial J}{\partial U}(x,U)$ is Lipschitz in $U$ on $\dom \times \cB_{2\gamma}(U_0)$ with Lipschitz constant $L_2$. Now set $\alpha \coloneq \frac{1}{L_2}$, and define  
    \begin{align}
        \begin{split} \label{eq:approx_convex_cost_def_Psi}
            \Psi \colon \dom \times \cB_{2\gamma}(U_0) & \to \Rm, \\ 
            (x,U) & \mapsto U - \alpha \frac{\partial J}{\partial U}(x,U). 
        \end{split}
    \end{align}
    Applying inequality \eqref{eq:gradient_descent_distance_x} from Lemma \ref{lem:gradient_descent} yields that for all $x \in \dom$ and $k \in \N$ 
    \begin{equation}\label{eq:approx_convex_cost_distPsiU} 
        \norm{\Psi^k(x, U_0) - U_0} \leq  \norm{\Psi^k(x, U_0) - U^*(x)} +  \norm{U^*(x) - U_0} \leq 2  \norm{U^*(x) - U_0} \leq  2 \gamma. 
    \end{equation}
    Thus, we have $\Psi^k(x, U_0) \in \cB_{2\gamma}(U_0)$ for all $k \in \N$ and the gradient step method given by $\Psi$ with starting point $U_0$ is well-defined on $\dom \times \cB_{2\gamma}(U_0)$. Further, by equation \eqref{eq:gradient_descent_distance_f} from \cref{lem:gradient_descent} for all $x \in \dom$ and $k \in \N$ it holds that 
    \begin{equation} \label{eq:approx_convex_cost_descent_J}
        J(x, \Psi^k(x, U_0)) - J(x, U^*(x)) \leq \frac{2 L_2}{k+4} \norm{U_0 - U^*(x)}^2 \leq \frac{2 L_2}{k+4} \gamma^2. 
    \end{equation} 
    Next, as in the proof of Theorem 6.1 in \cite{kang2022feedforward}, we use finite differences to approximate $\frac{\partial J}{\partial U}(x,U)$. For any $h >0$ define for $x \in \dom$, $U \in \cB_{2\gamma}(U_0)$, and $j \in \cu{1, \dots, m}$ 
    \begin{equation} \label{eq:approx_convex_cost_def_psi}
        \tilde{\psi}_{j,h}(x,U) \coloneq \frac{1}{h} \pr*{J(x, U_1, \dots, U_j+h, \dots, U_m) - J(x,U)}. 
    \end{equation}
    We have 
    \begin{equation*}
        \abs*{\frac{\partial J}{\partial U_j}(x,U) - \tilde{\psi}_{j,h}(x,U)} \leq C_j h, 
    \end{equation*}
    where $C_j \coloneq \max_{x \in \dom, U \in \cB_{2\gamma}(U_0)} \cu*{\frac{\partial^2 J}{\partial U_j^2}(x,U)}$. Hence, by defining 
    \begin{equation*}
        \tilde{\Psi}_h(x,U) \coloneq U - \alpha \begin{bmatrix}
            \tilde{\psi}_{1,h}(x,U), \dots, \tilde{\psi}_{m,h}(x,U)
        \end{bmatrix}^\top 
    \end{equation*}
    we obtain an approximation of $\Psi$ satisfying 
    \begin{equation} \label{eq:approx_convex_cost_error_tilde_psi}
        \norm{\Psi(x,U) - \tilde{\Psi}_h(x,U)}_2 \leq \alpha h \sqrt{m} \max_{j \in \cu{1, \dots, m}} C_j \leq \alpha h \sqrt{m} L_2 = h \sqrt{m}
    \end{equation}
    for all $(x,U) \in \dom \times \cB_{2 \gamma}(U_0)$, where we utilized the fact that for any matrix $X \in \R^{n \times n}$ it holds\footnote{For any $A \in \R^{n \times n}$ it is $\norm*{A} = \sup_{\norm*{y} = 1} \norm*{Ay} \geq \norm*{A e_j} =  \sqrt{\sum_{i = 1}^n A_{ij}^2} \geq \sqrt{A_{jj}^2} = \abs{A_{jj}}$.} that $\abs{X_{jj}} \leq \norm*{X}_2$ for all $j \in \cu{1, \dots, n}$. Let $\delta > 0$ and use the assumption on the existence of some NN $J_{\net, \delta}$ approximating $J$ to construct an approximation $\Psi^\net_{h, \delta}$ of $\tilde{\Psi}_h$ by substituting $J^\net_\delta$ for $J$ in the definition of $\tilde{\psi}_{j,h}$ in \eqref{eq:approx_convex_cost_def_psi} and then defining $\Psi^\net_{h, \delta}$ analogously to $\tilde{\Psi}_h$. Using \eqref{eq:approx_convex_cost_error_tilde_psi} and \eqref{eq:approx_convex_cost_NN_approx_J} we obtain 
    \begin{align}
        \begin{split} \label{eq:approx_convex_cost_NN_Psi}
            \norm*{\Psi(x,U) - \Psi^\net_{h, \delta}(x,U)} 
            & \leq \norm*{\Psi(x,U) - \tilde{\Psi}_h(x,U)} + \norm*{\tilde{\Psi}_h(x,U) - \Psi^\net_{h, \delta}(x,U)} \\
            & \leq \sqrt{m} h + 2 \frac{\alpha}{h} \sqrt{m} \delta. 
        \end{split}
    \end{align}
    Moreover, for all $(x,U) \in \dom \times \cB_{2\gamma}(U_0)$ it holds that 
    \begin{equation} 
        \norm*{I - \alpha \frac{\partial^2 J}{\partial U^2}(x,U)} \leq 1 
    \end{equation}
    since $\alpha = \frac{1}{L_2}$ and $\frac{\partial^2 J}{\partial U^2}(x,U) \succeq 0$, and thus each eigenvalue of $I - \alpha \frac{\partial^2 J}{\partial U^2}(x,U)$ is of the form $ 1 - \lambda$ for some $\lambda \in (0,1)$. Hence, we obtain that the mapping $\Psi$ defined in \eqref{eq:approx_convex_cost_def_Psi} is non-expansive on $\dom \times \cB_{2\gamma}(U_0)$. 
    Consequently, from \cref{lem:lipschitz_power_k} we get for all $(x,U) \in \dom \times \cB_{2\gamma}(U_0)$
    \begin{equation} \label{eq:approx_convex_cost_error_Psi_net}
        \norm*{\Psi^k(x,U) - \pr{\Psi^{\net}_{h, \delta}}^k(x,U)} \leq k \pr*{h \sqrt{m} + 2 \frac{\alpha}{h} \sqrt{m} \delta}. 
    \end{equation}
    Now let $k \in \N$, and let $h \in \R_{>0}, \delta \in \R_{>0}$ small enough such that \eqref{eq:approx_convex_cost_h_delta_small} is satisfied, i.e., $k \pr*{h \sqrt{m} + 2 \frac{\alpha}{h} \sqrt{m} \delta} \leq \gamma$. By \eqref{eq:approx_convex_cost_error_Psi_net} this yields together with \eqref{eq:approx_convex_cost_distPsiU} that for any $x \in \dom$
    \begin{equation} \label{eq:approx_convex_cost_distPsiNNU} 
        \norm*{U_0 - \pr{\Psi^{\net}_{h, \delta}}^k(x,U_0)} \leq \norm*{U_0 - \Psi^k(x,U_0)} + \norm*{\Psi^k(x,U_0) - \pr{\Psi^{\net}_{h, \delta}}^k(x,U)} \leq 3 \gamma.
    \end{equation}
    Now define the NN $U^{*\net}_{k,h,\delta}(x) \coloneqq \pr{\Psi^{\net}_{h, \delta}}^k(x, U_0)$ and observe that by \eqref{eq:approx_convex_cost_distPsiNNU} it holds that $U^{*\net}_{k,h,\delta}(\dom) \subset \cB_{3\gamma}(U_0) \subset \cU$. Applying \eqref{eq:approx_convex_cost_descent_J}, \eqref{eq:approx_convex_cost_error_Psi_net} and the fact that $J$ is Lipschitz in $U$ on $\dom \times \cB_{3\gamma}(U_0)$ with Lipschitz constant $L_1$ yields for all $x \in \dom$
    \begin{align*}
        J(x,U^{*\net}_{k,h,\delta}(x)) - J(x, U^*(x)) 
        & \leq \abs*{J(x,U^{*\net}_{k,h,\delta}(x) - J(x, \Psi^k(x))} + J(x,\Psi^k(x)) - J(x, U^*(x)) \\
        & \leq L_1 k \pr*{h \sqrt{m} + \frac{2}{h L_2} \sqrt{m} \delta} + \frac{2 L_2}{k+4} \gamma^2. 
    \end{align*}
    The number of neurons of the network $U_\net^*$ is given as $k$-times the number of neurons needed for the NN $\Psi_\net$, which equals $2 m \size \pr{J^{\net}_\delta}$. This concludes the proof. 
\end{mproof}

\subsection{Neural network approximations of optimal controls in non-strictly convex problems}
We now use Theorem~\ref{theo:approx_convex_cost} to state the main result of this section. The key step is to replace the assumption of the existence of a given neural network approximating $J$, as used in Theorem~\ref{theo:approx_convex_cost}, by explicitly constructing such an approximating network from the compositional structure of the problem data. Parts of the proof of the following theorem follow the proof of \cite[Theorem 6.2]{kang2022feedforward}.

\begin{theo}\label{theo:main_discete_convex}
    Consider a discrete-time optimal control problem as described in Setting \ref{set:general_dynamic} with $ \ell \equiv 0$. Assume that $\pr{f, \cG^f, \cL^f}$ and $\pr{g, \cG^g, \cL^g}$ are compositional functions that satisfy \cite[Assumption 1]{kang2022feedforward}. Further, assume that all the functions in the nodes of the compositional graph of $f$ and $g$ are twice continuously differentiable, respectively. Consider a compact set $\dom \subset \Rn$ of initial values and a convex set $\cU \subset \R^{qN}$ of admissible control values. Assume that the cost functional $J \colon \dom \times \cU \to \R, (x,U) \mapsto J(x,U)$ is  convex in $U$ for all $x \in \dom$. Further, suppose that there exists some $U_0 \in \cU$ and $\gamma > 0$ with $\cB_{3\gamma}(U_0) \subset \cU$ such that for all $x \in \dom$ there exists some $U^*(x) \in \cB_\gamma(U_0)$ satisfying 
    \begin{equation*}
        \inf_{U \in \cU} J(x, U) = J(x, U^*(x)). 
    \end{equation*} Then, for any $\varepsilon \in (0,1)$, there exists a NN $U^{*\net}_\varepsilon \colon \dom \to \cU$ that approximates $U^*$ in the weak sense, i.e., for all $x \in \dom$ the inequality
    \begin{equation*}
        J(x,U^{*\net}_\varepsilon(x)) - J(x, U^*(x)) \leq \varepsilon,
    \end{equation*}
    holds and the complexity of $U^{*\net}_\varepsilon$ is bounded by 
    \begin{equation*}
        \size \pr{U^{*\net}_\eps} \le C_1^r C_2 \, \pr*{\numberSteps \abs{ \cV_G^f } + \abs{ \cV_G^g}} \pr*{1 + \fC^r}  \frac{(qN)^{r+1}}{\varepsilon^{4r+1}}
    \end{equation*}
    where $C_1$ and $C_2$ are constants that depend at most polynomially on $\gamma, L_1,$ and $L_2$, and $\fC$ is a constant depending linearly on $L^g, C_f, L_{\max}^f, \Lambda^f, \abs{\cV_G^f}, C_g, L_{\max}^g, \Lambda^g, \abs{\cV_G^g}$, and polynomially of order $N$ on $L^f$.  
   Here $r \coloneqq \max \cu{r_{\max}^f, r_{\max}^g}$, $r_{\max}^f$ and $r_{\max}^g$ are part of the compositional features of the functions $f \colon \br{-R, R}^{n+q} \to \R^n$ and $g \colon \br{-R, R}^n \to \R$, $L^f$ and $L^g$ are Lipschitz constants of $f$ and $g$, %the functions $f \colon \br{-R, R}^{n+q} \to \R^n$ and $g \colon \br{-R, R}^n \to \R$, 
   respectively, and $R > 1$ is a constant large enough so that $\cB_{2 \gamma(U_0)} \subset [-R,R]^{qN}$ and 
   \begin{equation} \label{eq:main_discrete_convex:domain_f}
       f^k(x,\cdot) (U) \in \br*{-\frac{R}{2},\frac{R}{2}}^n
   \end{equation}
    for all $x \in \dom$, $U \in \cB_{2 \gamma}(U_0)$, and $0 \leq k \leq N-1$. Further, the Lipschitz constants $L_1$ and $L_2$ are defined as in Theorem \ref{theo:approx_convex_cost} and $C_f$ and $C_g$ are determined by $\cu{(d_{i,j}, m_{i,j}) \in \N^2 : f_{i,j} \in \cV_G^f}$ and $\cu{(d_{i,j}, m_{i,j}) \in \N^2 : g_{i,j} \in \cV_G^g}$, respectively.  
  \end{theo}

\begin{mproof}{\cref{theo:main_discete_convex}}
    Suppose that for each $\delta \in \R_{>0}$ we can prove the existence of a NN $J_\delta^\net \colon \Rn \times \R^{qN} \to \R$ with $\size(J_\delta^\net) \in \N$ neurons such that for all $(x,U) \in \dom \times \cB_{2 \gamma}(U_0)$ 
    \begin{equation} \label{eq:main_discrete_convex_NNs_approx_J}
        \lvert J(x, U) - J_{\delta}^\net(x,U) \rvert \leq \delta. 
    \end{equation}
     Now for $\varepsilon \in (0,1)$ we set 
    \begin{equation*}
        \bar{k} = \frac{6 L_2 \gamma^2}{\varepsilon}, 
        \quad 
        \bar{h} = \frac{1}{\bar{k} \sqrt{qN}} \min \cu*{\frac{\varepsilon}{3 L_1},\gamma},  
        \quad \text{and} \quad  
        \bar{\delta} =   \frac{\bar{h} L_2}{2 \sqrt{qN} \bar{k}} \min \cu*{\frac{\varepsilon}{3 L_1}, \gamma}. 
    \end{equation*}
    These choices of $\bar k, \bar h$ and $\bar \delta$ ensure that \eqref{eq:approx_convex_cost_h_delta_small} is satisfied. % and that the right-hand side in \eqref{eq:approx_convex_cost_main_inequality} is upper bounded by $\varepsilon$. 
    Thus, Theorem \ref{theo:approx_convex_cost} yields that there exists a NN $U^{*\net}_{\bar k, \bar h,\bar \delta} \colon \Rn \to \cU$ with $2 \bar{k} qN \size(J_{\bar \delta}^\net)$ neurons such that for all $x \in \dom$ we have 
    \begin{equation} \label{eq:main_discete_convex_JNN_inequality}
         J(x,U^{*\net}_{\bar k,\bar h,\bar \delta}(x)) - J(x, U^*(x)) \leq L_1 \bar{k} \pr*{\alpha \bar{h} \sqrt{qN} L_2 + 2 \frac{\alpha}{\bar{h}} \sqrt{qN} \bar{\delta}} + \frac{2 L_2}{\bar{k}+4} \gamma^2 \leq \varepsilon.  
    \end{equation}
   This yields the desired error bound. It is thus left to show the existence of the neural networks $J_\delta^\net$ \eqref{eq:main_discrete_convex_NNs_approx_J} and to estimate the total amount of needed neurons. Since $\pr{f, \cG^f, \cL^f}$ and $\pr{g, \cG^g, \cL^g}$ are compositional functions, we can construct neural networks as approximators of those functions by replacing each function corresponding to a node in the compositional function by a shallow NN with width $n_w \in \N$.
   By Theorem \ref{thm:NN_approximation_compositional}, for each $n_w \in \N$ there exist NNs $f^\net_{n_w}$ and $g^\net_{n_w}$ such that for all $(x,u) \in [-R,R]^{n+q}$  
   \begin{equation} \label{eq:main_discrete_convex:NN_f_g}
        \norm{f(x,u) - f^\net_{n_w}(x,u)}_2 \leq C_f L_{\max}^f \Lambda^f \abs{\cV_G^f} n_w^{ - \frac{1}{r_{\max}^f}} \quad \text{and} \quad \norm{g(x) - g^\net_{n_w}(x)}_2 \leq C_g L_{\max}^g \Lambda^g \abs{\cV_G^g} n_w^{ - \frac{1}{r_{\max}^g}},
   \end{equation}
   where $C_f$ and $C_g$ are determined by $\pr{d_{i,j}, m_{i,j} \mid f_{i,j} \in \cV_G^f}$ and $\pr{d_{i,j}, m_{i,j} \mid g_{i,j} \in \cV_G^g}$, respectively.
   Further, the number of neurons in $f^\net_{n_w}$ and $g^\net_{n_w}$ are given by $\size(f^\net_{n_w}) = \abs{\cV^f} n_w$ and $\size(g^\net_{n_w}) = \abs{\cV^g} n_w$, respectively. We now construct the NN $J_\delta^\net$ as composition of the nets $f^\net_{n_w}$ and $g^\net_{n_w}$ via  
   \begin{equation} \label{eq:main_discrete_convex:def_J_NN}
       J_\delta^\net (x,U) = g^\net_{n_w}(f^\net_{n_w}(x, \cdot)^\numberSteps (U)), \quad (x,U) \in \dom \times \cB_{2 \gamma}(U_0).  
   \end{equation}
   Supposing that 
   \begin{equation}  \label{eq:main_discrete_convex:domain_fNN}
        f_{n_w}^\net(x,\cdot)^k (U) \in \br*{-R, R}^n
   \end{equation}
   for all $(x,U) \in \dom \times \cB_{2 \gamma}(U_0)$ and $0 \leq k\leq N$, by \cite[Proposition 3.10]{kang2022feedforward} for all $n_w \in \N$ it holds for all $(x,U) \in \dom \times \cB_{2 \gamma}(U_0)$  
   \begin{align}
   \begin{split} \label{eq:main_discrete_convex:error_J}
        & \abs*{J(x,U) - g^\net_{n_w}( (f^\net_{n_w}(x,\cdot))^\numberSteps (U))} 
        =  \abs*{g( f^\numberSteps (x,\cdot) (U)) - g^\net_{n_w}( (f^\net_{n_w}(x,\cdot)^\numberSteps (U))} \\ 
       \leq &  L^g \frac{(L^f)^\numberSteps - 1 }{L^f - 1} C_f L_{\max}^f \Lambda^f \abs{\cV_G^f} n_w^{ - \frac{1}{r_{\max}^f}} 
       + C_g L_{\max}^g \Lambda^g \abs{\cV_G^g} n_w^{ - \frac{1}{r_{\max}^g}} 
       \leq \widetilde{C}_1 n_w^{ - \frac{1}{r}},
   \end{split}
   \end{align}
   where $r \coloneqq \max \cu{r_{\max}^f, r_{\max}^g}$ and $\widetilde{C}_1$ is a constant depending linearly on the constants
   \begin{equation*}
       L^g, C_f, L_{\max}^f, \Lambda^f, \abs{\cV_G^f}, C_g, L_{\max}^g, \Lambda^g, \abs{\cV_G^g}, 
   \end{equation*}
   as well as polynomially with order $\numberSteps$ on $L^f$. Now define 
   \begin{equation*}
       \fC \coloneqq \max\cu{C_1, 2 \frac{(L^f)^\numberSteps - 1 }{L^f - 1} C_f L_{\max}^f \Lambda^f \abs{\cV_G^f}}. 
   \end{equation*}
   and choose $n_w(\delta) \coloneqq \ceil*{\pr*{\nicefrac{\fC}{\delta}}^r}$. On the one hand, by recursively applying \cite[Proposition 3.10]{kang2022feedforward} and \eqref{eq:main_discrete_convex:NN_f_g} for $k = 0, \dots, N$ this choice yields for all $(x,U) \in \dom \times \cB_{2 \gamma}(U_0)$  
   \begin{equation*}
       \norm{f^k(x,\cdot) (U) - f^\net_{n_w(\delta)}(x,\cdot)^k (U)}_2 \leq \frac{1}{2} \leq \frac{R}{2},
   \end{equation*} 
   which together with  \eqref{eq:main_discrete_convex:domain_f} yields \eqref{eq:main_discrete_convex:domain_fNN}. Note that this iterative consideration is crucial to ensure that $ f^\net_{n_w(\delta)}(\cdot, \cdot)^k$ does not map to states $x$ outside of $[-R,R]^n$ where \eqref{eq:main_discrete_convex:NN_f_g} would no longer be applicable. On the other hand, from \eqref{eq:main_discrete_convex:error_J} we obtain the error bound for $J_\delta^\net$ proposed in \eqref{eq:main_discrete_convex_NNs_approx_J}. Thus, it remains to count the number of neurons in the network $ J_{\bar \delta}^\net$. It is given as 
   \begin{align*}
        \size(J_{\bar \delta}^\net) 
        = & \pr*{\numberSteps \abs{ \cV_G^f } + \abs{ \cV_G^g}} n_w(\bar \delta) \\
        \le &  \pr*{\numberSteps \abs{ \cV_G^f } + \abs{ \cV_G^g}} \pr*{ 1 + \pr[\Big]{\frac{\fC}{\bar \delta}}^{\!r}} \\
        \le & \pr*{ \numberSteps \abs{ \cV_G^f } + \abs{ \cV_G^g}} \pr*{1 + \fC^r} \bar{\delta}^{-r} \\
        = & \pr*{\numberSteps \abs{ \cV_G^f } + \abs{ \cV_G^g}} \pr*{1 + \fC^r} \pr*{\frac{1}{216} \frac{\varepsilon^3}{L_1 L_2 qN \gamma^4} \min \cu*{\frac{\varepsilon}{3 L_1}, \gamma}}^{-r} \\ 
        \leq & \pr*{\numberSteps \abs{ \cV_G^f } + \abs{ \cV_G^g}} \pr*{1 + \fC^r}  \frac{(qN)^r}{\varepsilon^{4r}} \pr*{\frac{1}{216}  \frac{1}{L_1 L_2 \gamma^4} \min \cu*{\frac{1}{3 L_1}, \gamma}}^{-r} \\
        =: & \pr*{\numberSteps \abs{ \cV_G^f } + \abs{ \cV_G^g}} \pr*{1 + \fC^r}  \frac{(qN)^r}{\varepsilon^{4r}} C_{1}^r. 
   \end{align*}
   In total, it follows that the number of neurons of the network $U^{*\net}_{\bar k, \bar h,\bar \delta}$ is given by 
   \begin{align*}
       \size(U^{*\net}_{\bar k, \bar h,\bar \delta}) 
       = 2 \bar{k} q N \size(J_{\bar \delta}^\net)
       \le & 2 \frac{6 L_2 \gamma^2}{\varepsilon} \, qN \, \pr*{\numberSteps \abs{ \cV_G^f } + \abs{ \cV_G^g}} \pr*{1 + \fC^r}  \frac{(qN)^r}{\varepsilon^{4r}} C_{1}^r \\
        = & C_1^r C_2 \, \pr*{\numberSteps \abs{ \cV_G^f } + \abs{ \cV_G^g}} \pr*{1 + \fC^r}  \frac{(qN)^{r+1}}{\varepsilon^{4r+1}}, 
   \end{align*}
   where $C_2 \coloneqq 12 L_2 \gamma^2$. This shows the claim.
\end{mproof}

In summary, Theorem~\ref{theo:main_discete_convex} establishes that, for an optimal control problem with compositional problem data and a convex cost functional, an approximate optimal feedback control can be realized using a number of neurons that does not grow exponentially with the state dimension provided that the constants appearing in Theorem~\ref{theo:main_discete_convex} also remain suitably bounded with the dimension. In addition to the Lipschitz constants of $g$, $f$, $J$, and $\frac{\partial J}{\partial U}$, this essentially requires that the growth of the compositional features of $g$ and $f$ remains suitably bounded with respect the dimension. Under this assumption, the curse of dimensionality can be mitigated. 
\begin{remark}
    Theorem \ref{theo:main_discete_convex} is formulated for discrete-time optimal control problems with terminal cost as we require $\ell \equiv 0$. An application to problems with stage costs can be achieved by using the reformulation described in Subsection \ref{subsec:ext_state_space_comp_struct} leveraging extend dynamics as in Remark \ref{rem:conv_reform_cost_funct}, given that this reformulation satisfies the assumptions of Theorem \ref{theo:main_discete_convex}.  
\end{remark}

\section{Conclusion} \label{sec:conclusion}
In this paper we have highlighted how compositional structures can be exploited in high-dimensional optimal control when using neural network approximations. By connecting neural network approximation theory for compositional functions with the analysis of control problems, we have identified conditions under which high-dimensional optimal controls can be represented without suffering from the curse of dimensionality. 
The developments presented here broaden the scope of existing results by relaxing convexity requirements and by showing how classical formulations of control problems can be adapted to fit those requirements. In this way, the theory is made applicable to a wider class of problems, including linear-quadratic settings and beyond.
For future research, a natural direction is to investigate how the theoretical representations derived here can inform practical algorithms for large-scale control. In particular, methods that automatically detect or exploit compositional features in nonlinear systems may provide a path toward efficient feedback synthesis with neural networks in high-dimensions.

\subsubsection*{Acknowledgements}
This work has been funded by the Deutsche Forschungsgemeinschaft (DFG, German \mbox{Research} Foundation) in the frame of the priority programme SPP 2298 ``Theoretical \mbox{Foundations} of Deep Learning'' – Project Number 464101154 and Project Number 463912816.

\bibliographystyle{acm}
\bibliography{bibcompfunc}

\end{document}

%% file: commands.tex
\DeclareFontEncoding{LS1}{}{}
\DeclareFontSubstitution{LS1}{stix}{m}{n}
\DeclareMathAlphabet{\mathscr}{LS1}{stixscr}{m}{n}

\newcommand*{\N}{\mathbb{N}}
\newcommand*{\R}{\mathbb{R}}

\newcommand*{\bI}{\mathbf{I}}

\newcommand{\cB}{\mathcal{B}}
\newcommand{\cC}{\mathcal{C}}

\newcommand{\cE}{\mathcal{E}}

\newcommand{\cG}{\mathcal{G}}

\newcommand{\cL}{\mathcal{L}}

\newcommand{\cU}{\mathcal{U}}
\newcommand{\cV}{\mathcal{V}}

\newcommand{\fC}{\mathfrak{C}}

\newcommand{\bff}{\mathbf{f}}
\newcommand{\bfg}{\mathbf{g}}

\newcommand{\cost}{\ell}
\newcommand{\terminal}{g}

\newcommand{\eps}{\varepsilon}

\newcommand{\diff}{D}
\newcommand{\dom}{\Omega}

\newcommand{\size}{\operatorname{size}}

\newcommand{\tnorm}[1]{{ \left\vert\kern-0.25ex \left\vert\kern-0.25ex \left\vert #1 \right\vert\kern-0.25ex \right\vert\kern-0.25ex \right\vert } }

\newcommand{\tnorms}[2]{ \mathopen{#1|\kern-0.25ex#1|\kern-0.25ex#1|}
#2
\mathclose{#1|\kern-0.25ex#1|\kern-0.25ex#1|}}

\DeclarePairedDelimiter{\cu}{\{}{\}}
\DeclarePairedDelimiter{\pr}{(}{)}
\DeclarePairedDelimiter{\br}{[}{]}
\DeclarePairedDelimiter{\abs}{\lvert}{\rvert}
\DeclarePairedDelimiter{\norm}{\lVert}{\rVert}

\DeclarePairedDelimiter{\ceil}{\lceil}{\rceil}

\newcommand{\SmallSum}[2]{ {\textstyle\sum\limits_{#1}^{#2}}}

\allowdisplaybreaks

%% file: bibcompfunc.bib
@article {kang2022feedforward,
    AUTHOR = {Kang, Wei and Gong, Qi},
     TITLE = {Feedforward neural networks and compositional functions with
              applications to dynamical systems},
   JOURNAL = {SIAM J. Control Optim.},
  FJOURNAL = {SIAM Journal on Control and Optimization},
    VOLUME = {60},
      YEAR = {2022},
    NUMBER = {2},
     PAGES = {786--813},
      ISSN = {0363-0129,1095-7138},
   MRCLASS = {41A30 (41A63 49M25 65P99 68T07)},
  MRNUMBER = {4395164},
MRREVIEWER = {Haijun\ Yu},
       DOI = {10.1137/21M1391596},
       URL = {https://doi.org/10.1137/21M1391596},
}

@book{nesterov2018lectures,
  title={Lectures on Convex Optimization},
  author={Nesterov, Yurii},
  volume={137},
  edition={2}, 
  year={2018},
  publisher={Springer}, 
  SERIES = {Springer Optimization and Its Applications},
  DOI = {10.1007/978-3-319-91578-4},
  URL = {https://doi.org/10.1007/978-3-319-91578-4},
}

@article{dahmen2025compositional_published,
	author = {Dahmen, Wolfgang},
	date = {2025/04/01},
	doi = {10.1007/s00365-024-09699-9},
	id = {Dahmen2025},
	isbn = {1432-0940},
    JOURNAL = {Constr. Approx.},
    FJOURNAL = {Constructive Approximation. An International Journal for
              Approximations and Expansions},
	number = {2},
	pages = {219--283},
	title = {Compositional Sparsity, Approximation Classes, and Parametric Transport Equations},
	url = {https://doi.org/10.1007/s00365-024-09699-9},
	volume = {61},
	year = {2025}}

@Article{poggio2017_WhyAndWhenCanDeepAvoidCoD,
  author    = {Tomaso Poggio and Hrushikesh Mhaskar and Lorenzo Rosasco and Brando Miranda and Qianli Liao},
  journal   = {International Journal of Automation and Computing},
  title     = {Why and when can deep-but not shallow-networks avoid the curse of dimensionality: A review},
  year      = {2017},
  issn      = {1751-8520},
  number    = {5},
  pages     = {503--519},
  volume    = {14},
  doi       = {10.1007/s11633-017-1054-2},
  publisher = {Springer Science and Business Media {LLC}},
}

@article {kang2024data,
    AUTHOR = {Kang, Wei and Sun, Kai and Xu, Liang},
     TITLE = {Data-driven computational methods for the domain of attraction
              and {Z}ubov's equation},
   JOURNAL = {IEEE Trans. Automat. Control},
  FJOURNAL = {Institute of Electrical and Electronics Engineers.
              Transactions on Automatic Control},
    VOLUME = {69},
      YEAR = {2024},
    NUMBER = {3},
     PAGES = {1600--1611},
      ISSN = {0018-9286,1558-2523},
   MRCLASS = {65L60 (93-08 93D20)},
  MRNUMBER = {4713397},
       DOI = {10.1109/tac.2023.3326769},
       URL = {https://doi.org/10.1109/tac.2023.3326769},
}

@Article{gong2023approximationReLU,
  author    = {Gong, Qi and Kang, Wei and Fahroo, Fariba},
  journal   = {Systems Control Lett.},
  FJOURNAL = {Systems \& Control Letters},
  title     = {Approximation of compositional functions with {ReLU} neural networks},
  year      = {2023},
  pages     = {105508},
  volume    = {175},
  publisher = {Elsevier},
}

@article{riekert2023deep,
  title={Deep neural network approximation of composite functions without the curse of dimensionality},
  author={Riekert, Adrian},
  journal={arXiv:2304.05790},
  eprint  = {arXiv:2304.05790},
  archivePrefix = {arXiv},
  year={2023}
}

@article{nakamura2021adaptive,
  author={Nakamura-Zimmerer, Tenavi and Gong, Qi and Kang, Wei},
  TITLE = {Adaptive deep learning for high-dimensional
              {H}amilton-{J}acobi-{B}ellman equations},
  JOURNAL = {SIAM J. Sci. Comput.},
  FJOURNAL = {SIAM Journal on Scientific Computing},
  volume={43},
  number={2},
  pages={A1221--A1247},
  year={2021},
  publisher={SIAM}
}

@article{nakamura2022,
  author={Nakamura-Zimmerer, Tenavi and Gong, Qi and Kang, Wei},
  journal={IEEE Open Journal of Control Systems}, 
  title={Neural Network Optimal Feedback Control With Guaranteed Local Stability}, 
  year={2022},
  volume={1},
  number={},
  pages={210-222},
  doi={10.1109/OJCSYS.2022.3205863}}

@article{Cybenko1989,
author="Cybenko, G.",
title="Approximation by superpositions of a sigmoidal function",
JOURNAL = {Math. Control Signals Systems},
  FJOURNAL = {Mathematics of Control, Signals, and Systems},
year="1989",
volume="2",
number="4",
pages="303--314",
issn="1435-568X",
doi="10.1007/BF02551274",
url="https://doi.org/10.1007/BF02551274"
}

@article{hornik1989multilayer,
  title={Multilayer feedforward networks are universal approximators},
  author={Hornik, Kurt and Stinchcombe, Maxwell and White, Halbert},
  journal={Neural networks},
  volume={2},
  number={5},
  pages={359--366},
  year={1989},
  publisher={Elsevier}
}

@incollection {Pinkus1999MLP,
    AUTHOR = {Pinkus, Allan},
     TITLE = {Approximation theory of the {MLP} model in neural networks},
 BOOKTITLE = {Acta numerica, 1999},
    SERIES = {Acta Numer.},
    VOLUME = {8},
     PAGES = {143--195},
 PUBLISHER = {Cambridge Univ. Press, Cambridge},
      YEAR = {1999},
      ISBN = {0-521-77088-2},
   MRCLASS = {41A30 (41A63 65D15 92B20)},
  MRNUMBER = {1819645},
MRREVIEWER = {Andrei\ Mart\'inez Finkelshtein},
       DOI = {10.1017/S0962492900002919},
       URL = {https://doi.org/10.1017/S0962492900002919},
}

@incollection {Guehring2023Expressivity,
    AUTHOR = {G\"uhring, Ingo and Raslan, Mones and Kutyniok, Gitta},
     TITLE = {Expressivity of deep neural networks},
 BOOKTITLE = {Mathematical aspects of deep learning},
     PAGES = {149--199},
 PUBLISHER = {Cambridge Univ. Press, Cambridge},
      YEAR = {2023},
      ISBN = {978-1-316-51678-2},
   MRCLASS = {68T07 (65Kxx)},
  MRNUMBER = {4505885},
}

@article{han2018solving,
  title={Solving high-dimensional partial differential equations using deep learning},
  author={Han, Jiequn and Jentzen, Arnulf and E, Weinan},
  JOURNAL = {Proc. Natl. Acad. Sci.},
  fjournal={Proceedings of the National Academy of Sciences},
  volume={115},
  number={34},
  pages={8505--8510},
  year={2018},
  publisher={National Academy of Sciences}
}

@Article{gonon2023deep,
  author    = {Gonon, Lukas and Schwab, Christoph},
  JOURNAL = {Anal. Appl. (Singap.)},
  FJOURNAL = {Analysis and Applications},
  title     = {Deep {ReLU} neural networks overcome the curse of dimensionality for partial integrodifferential equations},
  year      = {2023},
  number    = {01},
  pages     = {1--47},
  volume    = {21},
  publisher = {World Scientific},
}

@article{narendra1990identification,
  title={Identification and control of dynamical systems using neural networks},
  author={Narendra, Kumpati S and Parthasarathy, Kannan},
  journal={IEEE Trans. Neural Networks},
  fjournal={IEEE Transactions on neural networks},
  volume={1},
  number={1},
  pages={4--27},
  year={1990},
  publisher={IEEE}
}

@InProceedings{sontag1991feedback,
  author       = {Sontag, Eduardo D},
  booktitle    = {American Control Conference (ACC)},
  title        = {Feedback stabilization using two-hidden-layer nets},
  year         = {1991},
  organization = {IEEE},
  pages        = {815--820},
}

@article{kunisch2021semiglobal,
  title={Semiglobal optimal feedback stabilization of autonomous systems via deep neural network approximation},
  author={Kunisch, Karl and Walter, Daniel},
  JOURNAL = {ESAIM Control Optim. Calc. Var.},
  FJOURNAL = {ESAIM. Control, Optimisation and Calculus of Variations},
  volume={27},
  pages={16},
  year={2021},
  publisher={EDP Sciences},
  DOI = {10.1051/cocv/2021009},
  URL = {https://doi.org/10.1051/cocv/2021009},
}

@Article{gruene2021_ComputingLyapunovFunctions,
  author    = {Lars Gr{\"u}ne},
  journal   = {J. Comput. Dyn.},
  FJOURNAL = {Journal of Computational Dynamics},
  title     = {Computing {L}yapunov functions using deep neural networks},
  year      = {2021},
  number    = {2},
  pages     = {131--152},
  volume    = {8},
  publisher = {American Institute of Mathematical Sciences ({AIMS})},
  DOI = {10.3934/jcd.2021006},
  URL = {https://doi.org/10.3934/jcd.2021006},
}

@article{sperl2025_CLF,
title = {On the existence and neural network representation of separable control {L}yapunov functions},
journal = {Automatica},
volume = {182},
pages = {112517},
year = {2025},
issn = {0005-1098},
doi = {https://doi.org/10.1016/j.automatica.2025.112517},
url = {https://www.sciencedirect.com/science/article/pii/S0005109825004121},
author = {Mario Sperl and Jonas Mysliwitz and Lars Grüne},
}

@book {Bellman57,
    AUTHOR = {Bellman, Richard},
     TITLE = {Dynamic programming},
    SERIES = {Princeton Landmarks in Mathematics},
      NOTE = {Reprint of the 1957 edition,
              With a new introduction by Stuart Dreyfus},
 PUBLISHER = {Princeton University Press, Princeton, NJ},
      YEAR = {2010},
     PAGES = {xxx+340},
      ISBN = {978-0-691-14668-3},
   MRCLASS = {90-01 (49L20 90C39)},
  MRNUMBER = {2641641},
}

@incollection {NovakRitter1997,
    AUTHOR = {Novak, Erich and Ritter, Klaus},
     TITLE = {The curse of dimension and a universal method for numerical integration},
 BOOKTITLE = {Multivariate approximation and splines ({M}annheim, 1996)},
    SERIES = {Internat. Ser. Numer. Math.},
    VOLUME = {125},
     PAGES = {177--187},
 PUBLISHER = {Birkh\"{a}user, Basel},
      YEAR = {1997},
   MRCLASS = {65D30},
  MRNUMBER = {1485004},
}

@book {NovakWozniakowski2008I,
    AUTHOR = {Novak, Erich and Wo{\'z}niakowski, Henryk},
     TITLE = {Tractability of multivariate problems. {V}ol. 1: {L}inear
              information},
    SERIES = {EMS Tracts in Mathematics},
    VOLUME = {6},
 PUBLISHER = {European Mathematical Society (EMS), Z\"urich},
      YEAR = {2008},
     PAGES = {xii+384},
      ISBN = {978-3-03719-026-5},
   MRCLASS = {46E22 (41A63 41A65 46N40 47N40)},
  MRNUMBER = {2455266 (2009m:46037)},
MRREVIEWER = {Arthur G. Werschulz},
       DOI = {10.4171/026},
       URL = {http://dx.doi.org/10.4171/026},
}

@article{mnih2015human,
  title={Human-level control through deep reinforcement learning},
  author={Mnih, Volodymyr and Kavukcuoglu, Koray and Silver, David and Rusu, Andrei A and Veness, Joel and Bellemare, Marc G and Graves, Alex and Riedmiller, Martin and Fidjeland, Andreas K and Ostrovski, Georg and others},
  journal={Nature},
  volume={518},
  number={7540},
  pages={529--533},
  year={2015},
  publisher={Nature Publishing Group}
}

@article{silver2016mastering,
  title={Mastering the game of Go with deep neural networks and tree search},
  author={Silver, David and Huang, Aja and Maddison, Chris J and Guez, Arthur and Sifre, Laurent and Van Den Driessche, George and Schrittwieser, Julian and Antonoglou, Ioannis and Panneershelvam, Veda and Lanctot, Marc and others},
  journal={Nature},
  volume={529},
  number={7587},
  pages={484--489},
  year={2016},
  publisher={Nature Publishing Group}
}

@article{lillicrap2015continuous,
  title={Continuous control with deep reinforcement learning},
  author={Lillicrap, Timothy P and Hunt, Jonathan J and Pritzel, Alexander and Heess, Nicolas and Erez, Tom and Tassa, Yuval and Silver, David and Wierstra, Daan},
  journal={4th International Conference on Learning Representations, ICLR 2016},
  year={2016}
}

@article{hutzenthaler2019proof,
  author =        {Hutzenthaler, Martin and Jentzen, Arnulf and
                   Kruse, Thomas and Nguyen, Tuan Anh},
  fjournal =       {Partial Differential Equations and Applications},
  pages =         {1--34},
  publisher =     {Springer},
  JOURNAL = {Partial Differ. Equ. Appl.},
  title =         {A proof that rectified deep neural networks overcome
                   the curse of dimensionality in the numerical
                   approximation of semilinear heat equations},
  volume =        {1},
  year =          {2020},
}

@article{ChanMikaelWarin2019,
  author =        {Chan-Wai-Nam, Quentin and Mikael, Joseph and
                   Warin, Xavier},
  journal =       {J. Sci. Comput.},
  FJOURNAL = {Journal of Scientific Computing},
  number =        {3},
  pages =         {1667--1712},
  title =         {Machine learning for semi linear {PDE}s},
  volume =        {79},
  year =          {2019},
  doi =           {10.1007/s10915-019-00908-3},
  issn =          {0885-7474},
  url =           {https://doi.org/10.1007/s10915-019-00908-3},
}

@article{EHanJentzen2017CMStat,
  author =        {Weinan E and Jiequn Han and Arnulf Jentzen},
  journal =       {Commun. Math. Stat.},
  FJOURNAL = {Communications in Mathematics and Statistics},
  number =        {4},
  pages =         {349--380},
  title =         {Deep learning-based numerical methods for
                   high-dimensional parabolic partial differential
                   equations and backward stochastic differential
                   equations},
  volume =        {5},
  year =          {2017},
}

@article{BeckEJentzen2017,
  author =        {Christian Beck and Weinan E and Arnulf Jentzen},
  journal = {J. Nonlinear Sci.},
  fjournal =       {Journal of Nonlinear Science},
  number =        {4},
  pages =         {1563--1619},
  publisher =     {Springer},
  title =         {Machine learning approximation algorithms for
                   high-dimensional fully nonlinear partial differential
                   equations and second-order backward stochastic
                   differential equations},
  volume =        {29},
  year =          {2019},
  DOI = {10.1007/s00332-018-9525-3},
  URL = {https://doi.org/10.1007/s00332-018-9525-3},
}

@article {beck2021Kolmogorov,
    AUTHOR = {Beck, Christian and Becker, Sebastian and Grohs, Philipp and
              Jaafari, Nor and Jentzen, Arnulf},
     TITLE = {Solving the {K}olmogorov {PDE} by means of deep learning},
   JOURNAL = {J. Sci. Comput.},
  FJOURNAL = {Journal of Scientific Computing},
    VOLUME = {88},
      YEAR = {2021},
    NUMBER = {3},
     PAGES = {Paper No. 73, 28},
      ISSN = {0885-7474,1573-7691},
   MRCLASS = {65C30 (35R60 60H10 62M45 65M75)},
  MRNUMBER = {4293960},
       DOI = {10.1007/s10915-021-01590-0},
       URL = {https://doi.org/10.1007/s10915-021-01590-0},
}

@article {ReisingerZhang2020zerosum,
    AUTHOR = {Reisinger, Christoph and Zhang, Yufei},
     TITLE = {Rectified deep neural networks overcome the curse of
              dimensionality for nonsmooth value functions in zero-sum games
              of nonlinear stiff systems},
   JOURNAL = {Anal. Appl. (Singap.)},
  FJOURNAL = {Analysis and Applications},
    VOLUME = {18},
      YEAR = {2020},
    NUMBER = {6},
     PAGES = {951--999},
      ISSN = {0219-5305,1793-6861},
   MRCLASS = {65C30 (35R60 82C32)},
  MRNUMBER = {4154658},
MRREVIEWER = {Hayman\ Thabet},
       DOI = {10.1142/S0219530520500116},
       URL = {https://doi.org/10.1142/S0219530520500116},
}

@article{jentzen2025DNNBellman,
      title={Deep neural networks can provably solve {B}ellman equations for {M}arkov decision processes without the curse of dimensionality}, 
      author={Arnulf Jentzen and Konrad Kleinberg and Thomas Kruse},
      year={2025},
      journal={arXiv:2506.22851},
      eprint={2506.22851},
      archivePrefix={arXiv},
      primaryClass={math.OC},
      url={https://arxiv.org/abs/2506.22851}, 
}

@article {Gonon2024DeepStopping,
    AUTHOR = {Gonon, Lukas},
     TITLE = {Deep neural network expressivity for optimal stopping
              problems},
   JOURNAL = {Finance Stoch.},
  FJOURNAL = {Finance and Stochastics},
    VOLUME = {28},
      YEAR = {2024},
    NUMBER = {3},
     PAGES = {865--910},
      ISSN = {0949-2984,1432-1122},
   MRCLASS = {60G40 (60J05 62M45 65C30 68T07 91G20)},
  MRNUMBER = {4765842},
MRREVIEWER = {Scott\ Robertson},
       DOI = {10.1007/s00780-024-00538-0},
       URL = {https://doi.org/10.1007/s00780-024-00538-0},
}

@article{chen1994adaptively,
  title={Adaptively controlling nonlinear continuous-time systems using multilayer neural networks},
  author={Chen, Fu-Chuang and Liu, Chen-Chung},
  JOURNAL = {IEEE Trans. Automat. Control},
  FJOURNAL = {Institute of Electrical and Electronics Engineers.
              Transactions on Automatic Control},
  volume={39},
  number={6},
  pages={1306--1310},
  year={1994},
  publisher={IEEE}
}

@article{kunisch2023optimal,
  title={Optimal feedback control of dynamical systems via value-function approximation},
  author={Kunisch, Karl and Walter, Daniel},
  journal={Comptes Rendus. M{\'e}canique},
  volume={351},
  number={S1},
  pages={535--571},
  year={2023}
}

@article{sperl2025separable,
  title={Separable approximations of optimal value functions and their representation by neural networks},
  author={Sperl, Mario and Saluzzi, Luca and Kalise, Dante and Gr{\"u}ne, Lars},
  journal={arXiv:2502.08559},
  year={2025}
}

@article{Mhaskar1996NNoptimal,
    author = {Mhaskar, H. N.},
    title = {Neural Networks for Optimal Approximation of Smooth and Analytic Functions},
    journal = {Neural Computation},
    volume = {8},
    number = {1},
    pages = {164-177},
    year = {1996},
    month = {01},
    issn = {0899-7667},
    doi = {10.1162/neco.1996.8.1.164},
    url = {https://doi.org/10.1162/neco.1996.8.1.164},
    eprint = {https://direct.mit.edu/neco/article-pdf/8/1/164/813128/neco.1996.8.1.164.pdf},
}
